%% file: BoundedRatiosInFiniteAlgebras.tex
\documentclass{amsart}
\usepackage[utf8]{inputenc}

\title{Multiplicative Inequalities in Cluster Algebras of Finite Type}

\author{Michael Gekhtman}
\address{The University of Notre Dame\\
 255 Hurley Bldg\\ 46556 Notre Dame, Indiana, USA  }
\email{mgekhtma@nd.edu}

\author{Zachary Greenberg}
\address{Max Planck Institute for Mathematics in the Sciences\\
Inselstr. 22\\04103 Leipzig, 
Germany }
\email{greenberg@mis.mpg.de}

\author{Daniel Soskin}
\address{Institute for Advanced Study\\
1 Einstein Drive\\08540 Princeton, New Jersey, USA }
\email{dsoskin@ias.edu}

\usepackage{amsmath}
\usepackage{amsthm}
\usepackage{amssymb}
\usepackage{mathtools}
\usepackage{shuffle}
\usepackage{enumerate}
\usepackage{xifthen}
\usepackage{subcaption}
\usepackage[margin=1in]{geometry}

\usepackage{tikz}
\usepackage{tikz-cd}
 
\usepackage{dynkin-diagrams}
\pgfkeys{/Dynkin diagram,
edge length=.25cm,}

\usetikzlibrary {arrows.meta}
\usetikzlibrary{decorations.pathreplacing,decorations.markings}
\usetikzlibrary{shapes.geometric,quotes}
\usepackage{pgfplots}
\usepgfplotslibrary{fillbetween}
\pgfplotsset{compat=1.9}

\usepackage{hyperref}
\hypersetup{
	colorlinks=true,
	linkcolor=blue,
	hypertexnames=false,
	filecolor=black,      
	urlcolor=blue,
	citecolor=blue
}
\usepackage[nameinlink]{cleveref}
\usepackage{blkarray}
\usepackage{nicematrix}
\usepackage{lscape}
\usepackage[title]{appendix} 

\usetikzlibrary{positioning, calc}
\usetikzlibrary{arrows.meta}
\usetikzlibrary{decorations.markings}
\usetikzlibrary{shapes.geometric}
\tikzstyle{mutable}=[circle,draw=black,fill=black,inner sep=0pt,minimum size=6]
\tikzstyle{frozen}=[draw=black,fill=blue,inner sep=0pt,minimum size=6]
\tikzstyle{mutableBig}=[circle,draw=black,double=white,fill=black,inner sep=0pt,minimum size=10]
\tikzstyle{frozenBig}=[draw=black,fill=black,inner sep=0pt,minimum size=10]
\tikzstyle{barrow}=[->,very thick,>=Triangle]
\newcommand{\distL}{1.2}

\newcommand{\st}{\mid}
\newcommand{\keyword}[1]{\textbf{\emph{#1}}}
\newtheorem{theorem}{Theorem}[section]
\newtheorem{lemma}[theorem]{Lemma} 
\newtheorem{claim}[theorem]{Claim}
\newtheorem{proposition}[theorem]{Proposition}
\newtheorem{corollary}[theorem]{Corollary}
\newtheorem{conjecture}[theorem]{Conjecture}
\theoremstyle{definition}
\newtheorem{remark}[theorem]{Remark}
\newtheorem{definition}[theorem]{Definition}
\newtheorem{example}[theorem]{Example}
\newtheorem{problem}[theorem]{Problem}

\newcommand{\sslash}{\mathbin{/\mkern-6mu/}}

\newcommand{\Z}{\mathbb{Z}}

\newcommand{\R}{\mathbb{R}}

\newcommand{\C}{\mathbb{C}}

\newcommand{\spn}{\mathrm{span}}

\newcommand{\tensor}{\otimes}
\renewcommand{\vec}[1]{\mathbf{#1}}
\DeclareMathOperator{\Hom}{Hom}

\DeclareMathOperator{\rank}{Rank}
\newcommand{\cone}[1]{#1}

\newcommand{\semifield}[1]{#1}
\newcommand{\multgroup}[1]{#1^\times}
\newcommand{\ideal}[1]{\langle #1 \rangle}

\newcommand{\exchangeMat}[1]{#1}
\newcommand{\exchangeMatFull}[1]{\widehat{#1}}
\newcommand{\exchangeMatExtended}[1]{\widetilde{#1}}
\newcommand{\plus}[1]{[#1]_{+}}
\newcommand{\minus}[1]{[#1]_{-}}

\newcommand{\xvar}[1]{#1}
\newcommand{\xindexSet}{\Pi}
\newcommand{\yvar}[1]{#1}
\newcommand{\uvar}[1]{#1}

\newcommand{\clusterAlg}[1]{\mathcal{#1}}
\newcommand{\A}{\clusterAlg{A}}
\newcommand{\localizedClusterAlg}[1]{\clusterAlg{#1}[\xvar{\vec{x}^{-1}}]}
\newcommand{\pospoints}{X_{>0}}

\newcommand{\yAlg}[1]{\mathcal{#1}}

\newcommand{\uAlg}[1]{\mathcal{#1}}
\newcommand{\uPoints}{\mathcal{M}_D}
\newcommand{\uPosPoints}{\mathcal{M}_D^{>0}}

\newcommand{\weight}{\mathrm{wt}}
\newcommand{\ratiovec}[1]{\mathbf{#1}}

\renewcommand{\root}[1]{#1}
\newcommand{\compDegree}[2]{\epsilon(#1,#2)}

\newcommand{\xgamma}{\xvar{x}_{\root{\gamma}}}
\newcommand{\xomega}{\xvar{x}_{\root{\omega}}}

\newcommand{\gr}{\mathrm{Gr}}
\newcommand{\pl}[1]{[#1]}

\begin{document}

\begin{abstract}
    Generalizing the notion of a multiplicative inequality among minors of a totally positive matrix, we describe, over full rank cluster algebras of finite type, the cone of Laurent monomials in cluster variables that are bounded as a real-valued function on the  positive locus of the cluster variety. 
    We prove that the extreme rays of this cone are the $u$-variables of the cluster algebra. Using this description, we prove that all bounded ratios are bounded by 1 and give a sufficient condition for all such ratios to be subtraction free. This allows us to show in  $\gr(2,n)$, $\gr(3,6), \gr(3,7), \gr(3,8)$ that every bounded Laurent monomial in Pl\"ucker coordinates factors into a positive integer combination of so-called primitive ratios. In $\gr(4,8)$ this factorization does not exists, but we provide the full list of extreme rays of the cone of bounded Laurent monomials in Pl\"ucker coordinates. 
\end{abstract}

\maketitle

\section{Introduction}

Totally nonnegative matrices are matrices in which each minor is nonnegative. These matrices appear in a wide range of mathematical areas including higher Teichm\"uller theory and the representation theory of quantum groups. In fact the search for minimal sets of minors to guarantee a matrix is totally nonnegative was one of the problems that inspired the theory of cluster algebras \cite{fz-doubleBruhat}. This search was also related to describing the {\em dual canonical bases} in the representation theory of quantum groups. Lusztig has shown that specialization of elements of the dual canonical basis at q=1 are totally nonnegative polynomials, which are polynomials in matrix entries attaining nonnegative values on totally nonnegative matrices \cite{LusztigTPCB}. To this end, there is an interest in determinantal inequalities which are a natural source of such polynomials.
\begin{definition}
    \keyword{Determinantal inequalities} are inequalities in real linear combinations of products of minors which hold over all totally nonnegative matrices.
\end{definition} 

Determinantal inequalities have been studied for years, starting with classical results by Hadamard, Fischer and Koteljanskii  \cite{hadamard1893}, \cite{fischer1908}, \cite{KotelTheory}, \cite{KotelTheoryRuss}. A determinantal inequality is called multiplicative when it compares two products of minors. The problem of describing the set of multiplicative determinantal inequalities for totally positive matrices was stated by S. Fallat and C. Johnson in \cite{fallat2000det} in terms of bounded ratios of products of minors. 

\begin{problem}\label{p:R} Describe ratios $R$ of products of minors bounded as a real-valued function on the locus of totally positive elements in $GL_{n}$, where $R$ is of the form
\begin{equation}\label{eq:R}{\det(A}_{I_{1},I'_{1}}){\det(A}_{I_{2},I'_{2}})...{\det(A}_{I_{p},I'_{p}})/{\det(A}_{J_{1},J'_{1}}) {\det(A}_{J_{2},J'_{2}})...{\det(A}_{J_{q},J'_{q}})
\end{equation}   
where $I_k,I'_k,J_k,J'_k \subseteq \{1,2,\ldots,n\}$ with
$|I_k|=|I'_k|$ and $|J_k|=|J'_k|$.
\end{problem} 

In \cite{fallat2003multiplicative},  necessary
and sufficient conditions were given for a ratio of products of two principal minors to
be bounded over totally positive matrices. 

This result was generalized to non-principal minors in
\cite{skandera2004inequalities} by M. Skandera, whose approach then allowed A. Boocher and B. Froehle to restate the problem in terms of ratios of products of Pl\"ucker coordinates bounded over the totally positive Grassmannian \cite{boocher2008generators}. 
The necessary condition in
\cite{fallat2003multiplicative} was generalized in \cite{boocher2008generators} to the case of non-principal minors and an explicit factorization of ratios of products of two minors
into products of so-called {\em primitive ratios} was constructed.
It has been conjectured that all bounded ratios can be factored into products of nonnegative integer powers 
of primitive ratios
\cite{fallat2003multiplicative}.
Recently, it was proved in \cite{soskin2023bounded} that for any $n\times n$ matrix the set of bounded ratios is finitely generated. Along with that some examples of nonprimitive generators were discovered which disprove the conjecture on factorization into primitive ratios for matrices of order $n\geq 4$ stated in \cite{fallat2003multiplicative}. 
Other conjectures on bounded ratios stated in \cite{boocher2008generators} remain open. 
\begin{conjecture}\label{c:1} 
 Ratio (\ref{eq:R}) is bounded if and only if it is bounded by 1.  
\end{conjecture}
Both totally positive matrices and elements of totally positive Grassmannian can be parameterized using weighted planar networks \cite{fz-doubleBruhat, postnikov2006total}. In this parametrization, any minor (resp. Pl\"ucker coordinate) is a polynomial in terms of {\em face weights}.
\begin{definition}\label{df:subfreeFW}
A ratio $\frac{p}{q}$ is called \emph{subtraction-free (in positive weights)} if $q-p$ is a polynomial function in face weights with all positive coefficients.
\end{definition}
\begin{conjecture}\label{c:2} Ratio (\ref{eq:R}) is bounded if and only if it is subtraction free.
\end{conjecture} 

Recall that every Pl\"ucker coordinate belongs to the set of cluster variables in the cluster structure on a Grassmannian \cite{scott-grassmannians, gekhtman2010cluster}. On the other hand, the totally positive locus is defined for any cluster algebra of geometric type. Thus, it is natural to extend \Cref{p:R} to ratios in all cluster variables in a given cluster algebra. 
\begin{problem}\label{p:clusterR} Describe ratios of products of cluster variables bounded as a real-valued function on the totally positive locus.
\end{problem} 

For brevity, we will simply call a ratio of products of cluster variables (equivalently, a Laurent monomial in cluster variables) which is bounded on the totally positive locus {\em a bounded ratio}.)

In this note we characterize the cone of bounded ratios in \Cref{p:clusterR} for full rank cluster algebras of finite type. 
\begin{theorem}
    Let $\clusterAlg{D}$ be a full rank finite type cluster algebra associated to a Dynkin diagram $D$. The generators of the cone of bounded ratios are of the form \[ \frac{\prod\limits_{\root{\gamma} \rightarrow \root{\omega} } \xvar{x}_{\root{\omega}}}{\xgamma \xgamma'}\] where $\xgamma$ is the variable at the source of a Dynkin type bipartite quiver that mutates to $\xgamma'$.  
\end{theorem}
\begin{corollary}
    The generators of the bounded cone correspond exactly to the u-variables of the cluster algebra (see \cite{ahsl-ClusterConfigurationOfFiniteType} for definition). 
\end{corollary}

This classification allows us to prove the analogue of \Cref{c:1} for full rank cluster algebras of finite type. 
\begin{corollary}
    Every bounded ratio in cluster variables is bounded by $1$.
\end{corollary}

The analogue of \Cref{df:subfreeFW} is 
\begin{definition}
    A ratio $\frac{p}{q}$ is \keyword{subtraction free (in cluster variables)} if $q-p$ can be expressed as a polynomial function with positive coefficients in the full set of cluster variables. Note that this implies that when $q-p$ is expressed as a Laurent polynomial in the cluster variables from a fixed cluster, the coefficients are all positive. 
\end{definition}
We provide an explicit counter example showing the following natural analogue of \Cref{c:2} is false in general (\Cref{ex:BoundedButNotSubtractionFree}):
\begin{conjecture}
    Every bounded ratio is subtraction free in cluster variables.
\end{conjecture}
However we show the following:
\begin{corollary}
    Every \textbf{extreme ray} of the bounded cone is subtraction free when expressed as a polynomial in all cluster variables.
\end{corollary} 
This proves that every bounded ratio that can be factored into integer powers of extreme rays is subtraction free in cluster variables. Thus the only counter examples are when the integer powers of extreme rays miss some bounded integer ratios. 

In type $A_{n-3}$ every cluster variable can be interpreted as a Pl\"ucker coordinate  in $\gr(2,n)$. However in other types there are other kinds of cluster variables. Thus we provide an algorithm to reduce the full bounded cone to the bounded cone generated by any finite subset of the cluster variables. Using this we study the cone of Pl\"ucker bounded ratios in $\gr(3,6)$, $\gr(3,7)$, $\gr(3,8)$. We show 
\begin{theorem}
    Every bounded ratio of Pl\"ucker coordinates in  $\gr(2,n),\gr(3,6), \gr(3,7)$ and $\gr(3,8)$ is bounded by 1 and is subtraction free in face weights. Furthermore every such ratio can be factored into a positive integer combination of primitive ratios. 
\end{theorem}
As discussed in \cite{soskin2023bounded} not every Pl\"ucker bounded ratio can be factored into a positive combination of primitive ratios. In \Cref{fig:Gr48ExtremeRays} we provide the full list of extreme rays of the cone. We checked that every such ratio is subtraction free in positive weights and thus every bounded ratio in Pl\"ucker coordinates on $\gr(4,8)$ is bounded by 1. It remains open in this case if every such ratio can be factored as a positive \emph{integer} combination of these extreme rays.

The paper proceeds as follows. In \Cref{sec:ClusterAlgebras} we review basic facts from cluster algebra theory. 
In \Cref{sec:BoundedRatios} we introduce bounded ratios in cluster variables and torus action which suggests a necessary conditions for a ratio to be bounded. Furthermore, we show that generators of the cone of bounded ratios in cluster variables are in correspondence with u-variables. In \Cref{sec:SubsetRatios}, we discuss bounded ratios in any fixed subset of cluster variables along with examples for $\gr(3, 6)$, $\gr(3,7)$ and $\gr(3,8)$.  
 
\section{Cluster Algebras}\label{sec:ClusterAlgebras}
In this section we review the definitions of a cluster algebra. For a complete introduction see \cite{fwz-IntroChap1-3}. 

\begin{definition}
    A square matrix $B$ is \keyword{skew symmetrizable} if there is a diagonal matrix $W$ with positive integer entries such that $WB$ is skew symmetric. 
\end{definition}
\begin{definition}
    A \keyword{seed} of a cluster algebra consists of a pair $(\exchangeMatFull{B},\xvar{\vec{x}})$ where $\exchangeMatFull{B}$ is a $(n+m) \times (n+m)$ skew symmetrizable matrix called the \keyword{full exchange matrix} and $\xvar{\vec{x}} = \{\xvar{x}_1,\dots,\xvar{x}_n,\xvar{f}_1,\dots,\xvar{f}_m\}$ is a collection of $n+m$ commuting variables called \keyword{cluster variables}. The first $n$ variables $\xvar{x}_1,\dots, \xvar{x}_n$ are called \keyword{mutable} or \keyword{unfrozen} and the last $m$ variables $\xvar{f}_1,\dots,\xvar{f}_m$ are called \keyword{frozen}.
\end{definition}
\begin{remark}
     We consider the rows and columns of $\exchangeMatFull{B}$ to be indexed by the cluster variables. Thus we call an index $i$ frozen or unfrozen if the corresponding cluster variable is frozen or unfrozen respectively. 
\end{remark}

\begin{remark}
    It is convenient to represent $\exchangeMatFull{B}$ with a quiver $Q$. When $\exchangeMatFull{B}$ is skew symmetric, $Q$ has a node for each row/column and an arrow of weight $\exchangeMatFull{B}_{ij}$ from $i$ to $j$ when $\exchangeMatFull{B}_{i,j}>0$. When $\exchangeMatFull{B}$ is skew symmetrizable with matrix $W$ the associated quiver has weighted nodes. See Section 3 of \cite{zickert-FockGoncharovForRankTwoLieGroups} for full details. We will often borrow terminology of quivers to refer to exchange matrices. For example, we refer to indices of $\exchangeMatFull{B}$ as \keyword{sources/sinks} if the node of the associated quiver is a source or sink of the subquiver of mutable nodes. 
\end{remark}

\begin{definition}
    The \keyword{mutable part} of a quiver/matrix is the subgraph/submatrix indexed by mutable variables. This corresponds to a matrix $\exchangeMat{B}$ consisting of the top left $n\times n$ block of $\exchangeMatFull{B}$. We also consider the \keyword{extended exchange matrix}, $\exchangeMatExtended{B}$, which is the rectangular $n\times (n+m)$ submatrix of the mutable rows of $\exchangeMatFull{B}$.
\end{definition}

\begin{definition}
    A \keyword{$Y$-seed}, $(\exchangeMatFull{B},\vec{\yvar{y}})$, is a pair of skew symmetrizable matrix and list of variables $\vec{\yvar{y}} = \{\yvar{y}_1,\dots,\yvar{y}_n\}$, called $Y-$variables. Note that a $Y-$seed only has a variable for each mutable index and thus only depends on $B$.
\end{definition}

\begin{definition}
    We say a seed is \keyword{full rank} if $\rank(\exchangeMatExtended{B}) = n$. 
\end{definition}
\begin{remark}
    We note that the property of being full rank depends on the frozen nodes. In fact every seed can be made full rank by ``framing'' each mutable node with corresponding frozen node attached out from the mutable node.
\end{remark}

\begin{example}
    In \Cref{fig:A1seeds} we see a variety of seeds represented as quivers and exchange matrices. The two seeds in \Cref{fig:A1seedNoFrozen} and \Cref{fig:A1seedFullRank} have the same mutable part. As such we refer to these seeds as type $A_1$. Although they are the same type the seed in \Cref{fig:A1seedNoFrozen} is not full rank, while the seed in \Cref{fig:A1seedFullRank} is. In \Cref{fig:B3seedNoFrozen} and \Cref{fig:C3seedNoFrozen} we see two examples with skew symmetrizable matrices. In quivers we use a large node with an extra ring to represent nodes of weight 2.
    \begin{figure}[hb]
        \centering
        \begin{subfigure}[]{.24\textwidth}
            \centering
            \begin{equation*}
               \begin{blockarray}{c c}
                   & \xvar{x}_1\\
                    \begin{block}{c [c ]}
                      \xvar{x}_1 &  0 \\
                    \end{block}\\
                \end{blockarray}                             
            \end{equation*}
            \begin{tikzpicture}
                \node[mutable]	(X1) at (0,0) [label = above: $\xvar{x}_1$]	{};
            \end{tikzpicture}
        \caption{$A_1$}
        \label{fig:A1seedNoFrozen}
        \end{subfigure}
        \begin{subfigure}[]{.24\textwidth}
            \centering
            \begin{equation*}
               \begin{blockarray}{c c  c}
                   & \xvar{x}_1 & \xvar{f}_1\\
                    \begin{block}{c [c  c]}
                      \xvar{x}_1 &  0  & 1 \\
                     \xvar{f}_1 & -1 & 0\\
                    \end{block}
                \end{blockarray}                             
            \end{equation*}
            \begin{tikzpicture}
                \node[mutable]	(X1) at (0,0) [label = above: $\xvar{x}_1$]	{};
                \node[frozen]	(F1) at (\distL,0) [label = above: $\xvar{f}_1$] {};
        		\draw[barrow] (X1) to (F1);
            \end{tikzpicture}
        \caption{$A_1$}
        \label{fig:A1seedFullRank}
        \end{subfigure}
        \begin{subfigure}[]{.24\textwidth}
            \centering
               \begin{equation*}
               \begin{blockarray}{c c  c }
                    \xvar{x}_1 & \xvar{x}_2 & \xvar{x}_3\\
                    \begin{block}{[c c  c]}
                       0 & 1 &0\\
                      -1 & 0 & -2 \\
                       0 & 1 & 0\\
                    \end{block}
                \end{blockarray}                             
            \end{equation*}
            \begin{tikzpicture}
                \node[mutable]	(X1) at (0,0) [label = above: $\xvar{x}_1$]	{};
                \node[mutable]	(X2) at (\distL,0) [label = above: $\xvar{x}_2$] {};
                \node[mutableBig]	(X3) at (2*\distL,0) [label = above: $\xvar{x}_3$]	{};
                \draw[barrow] (X1) to (X2);
                \draw[barrow] (X3) to (X2);
            \end{tikzpicture}
        \caption{$B_3$}
        \label{fig:B3seedNoFrozen}
        \end{subfigure}
        \begin{subfigure}[]{.24\textwidth}
            \centering
               \begin{equation*}
               \begin{blockarray}{c c  c }
                    \xvar{x}_1 & \xvar{x}_2 & \xvar{x}_3\\
                    \begin{block}{[c c  c]}
                       0 & 1 &0\\
                      -1 & 0 & -1 \\
                       0 & 2 & 0\\
                    \end{block}
                \end{blockarray}                             
            \end{equation*}
            \begin{tikzpicture}
                \node[mutableBig]	(X1) at (0,0) [label = above: $\xvar{x}_1$]	{};
                \node[mutableBig]	(X2) at (\distL,0) [label = above: $\xvar{x}_2$] {};
                \node[mutable]	(X3) at (2*\distL,0) [label = above: $\xvar{x}_3$]	{};
                \draw[barrow] (X1) to (X2);
                \draw[barrow] (X3) to (X2);
            \end{tikzpicture}
        \caption{$C_3$}
        \label{fig:C3seedNoFrozen}
        \end{subfigure}
        \caption{Seeds as Matrices and Quivers}
        \label{fig:A1seeds}
    \end{figure}
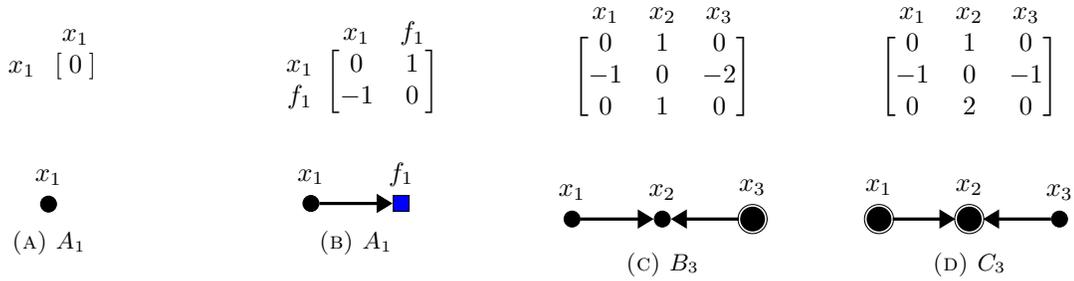
\end{example}

\begin{example}
    Our key running example will be cluster algebras of type $A_3$. In \Cref{fig:A3seedNoFrozen} we see the version with only mutable nodes. This seed is also not full rank. To make a full rank seed with the same mutable part it suffices to add one frozen node attached to either $x_1$ or $x_3$ breaking the symmetry. However it will be convenient to add more frozen nodes to obtain a quiver for cluster algebra associated to $\gr(2,6)$ as in \Cref{fig:A3seedGr26}.
    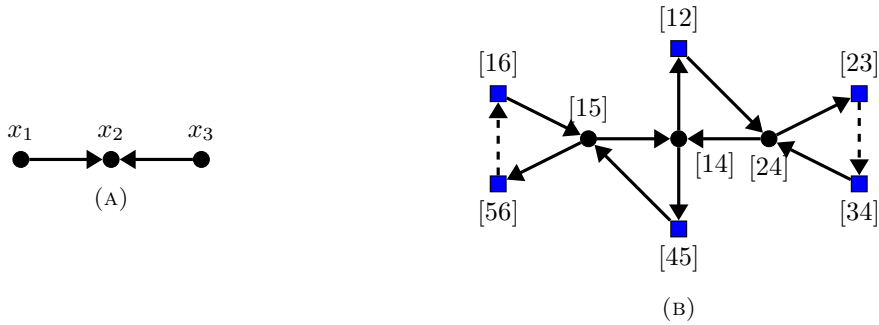
\begin{figure}[hb]
        \centering
        \begin{subfigure}[]{.3\textwidth}
            \centering
            \begin{tikzpicture}
                \node[mutable]	(X1) at (0,0) [label = above: $\xvar{x}_1$]	{};
                \node[mutable]	(X2) at (\distL,0) [label = above: $\xvar{x}_2$] {};
                \node[mutable]	(X3) at (2*\distL,0) [label = above: $\xvar{x}_3$]	{};
                \draw[barrow] (X1) to (X2);
                \draw[barrow] (X3) to (X2);
            \end{tikzpicture}
        \caption{}
        \label{fig:A3seedNoFrozen}
        \end{subfigure}
        \begin{subfigure}[]{.6\textwidth}
            \centering
            \begin{tikzpicture}
                \node[mutable]	(X15) at (2*\distL-\distL,0) [label = above: $\pl{15}$]	{};
                \node[mutable]	(X14) at (2*\distL+0,0) [label = below right: $\pl{14}$]	{};
                \node[mutable]	(X24) at (2*\distL+\distL,0) [label = below: $\pl{24}$]	{};
                \node[frozen]	(F12) at (2*\distL+0,\distL) [label = above: $\pl{12}$] {};
                \node[frozen]	(F23) at (2*\distL+2*\distL,0.5*\distL) [label = above: $\pl{23}$] {};
                \node[frozen]	(F34) at (2*\distL+2*\distL,-0.5*\distL) [label = below: $\pl{34}$] {};
                \node[frozen]	(F45) at (2*\distL+0,-\distL) [label = below: $\pl{45}$] {};
                \node[frozen]	(F56) at (2*\distL-2*\distL,-0.5*\distL) [label = below: $\pl{56}$] {};
                \node[frozen]	(F16) at (2*\distL-2*\distL,0.5*\distL) [label = above: $\pl{16}$] {};
        		\draw[barrow] (X15) to (X14);
                \draw[barrow] (X24) to (X14);
                \draw[barrow] (X14) to (F12);
                \draw[barrow] (F12) to (X24);
                \draw[barrow] (X24) to (F23);
                \draw[barrow] (F34) to (X24);
                \draw[barrow] (X14) to (F45);
                \draw[barrow] (F45) to (X15);
                \draw[barrow] (X15) to (F56);
                \draw[barrow] (F16) to (X15);
                \draw[barrow,dashed] (F56) to (F16);
                \draw[barrow,dashed] (F23) to (F34);
            \end{tikzpicture}
        \caption{}
        \label{fig:A3seedGr26}
        \end{subfigure}
        \caption{Seeds of Type $A_3$}
        \label{fig:A3seeds}
    \end{figure}
\end{example}

\begin{figure}[hb]
    \centering
    {\renewcommand{\distL}{1.8}
     \begin{tikzpicture}
                \node[frozen]	(F456) at (0*\distL,0*\distL) [label = above: $\pl{456}$]	{};
                \node[mutable]	(X356) at (1*\distL,0*\distL) [label = above: $\pl{356}$]	{};
                \node[mutable]	(X256) at (2*\distL,0*\distL) [label = above: $\pl{256}$]	{};
                \node[frozen]	(F156) at (3*\distL,0*\distL) [label = above: $\pl{156}$]	{};
                \node[mutable]	(X346) at (1*\distL,-1*0.75*\distL) [label = left: $\pl{346}$]	{};
                \node[mutable]	(X236) at (2*\distL,-1*0.75*\distL) [label = below left: $\pl{236}$]	{};
                \node[frozen]	(F126) at (3*\distL,-1*0.75*\distL) [label = below left: $\pl{126}$]	{};
                \node[frozen]	(F345) at (1*\distL,-2*0.75*\distL) [label = below: $\pl{345}$]	{};
                \node[frozen]	(F234) at (2*\distL,-2*0.75*\distL) [label = below: $\pl{234}$]	{};
                \node[frozen]	(F123) at (3*\distL,-2*0.75*\distL) [label = below: $\pl{123}$]	{};
                \draw[barrow] (F456) to (X356);
                \draw[barrow] (X356) to (X256);
                \draw[barrow] (X256) to (F156);
                \draw[barrow] (X346) to (X236);
                \draw[barrow] (X236) to (F126);
                \draw[barrow] (X356) to (X346);
                \draw[barrow] (X346) to (F345);
                \draw[barrow] (X256) to (X236);
                \draw[barrow] (X236) to (F234);
                \draw[barrow,dashed] (F156) to (F126);
                \draw[barrow,dashed] (F126) to (F123);
                \draw[barrow] (X236) to (X356);
                \draw[barrow] (F126) to (X256);
                \draw[barrow] (F234) to (X346);
                \draw[barrow] (F123) to (X236);
                \draw[barrow, bend left,dashed] (F345) to (F456);
            \end{tikzpicture}
            }
    \caption{Grassmannian Seed for $Gr(3,6)$}
    \label{fig:GrassmanSeed}
\end{figure}
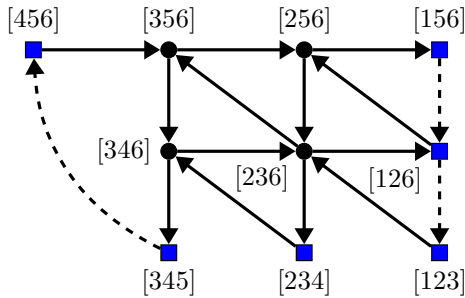

\begin{example}
    Our other key example will be seeds corresponding to the Grassmannian $\gr(k,n)$. In \cite{fwz-IntroChap6} Section 6.7 they give a construction of an initial seed with $(k-1)(n-k-1)$ mutable nodes and $n$ frozen nodes arranged in an $k \times (n-k)$ grid. In \Cref{fig:GrassmanSeed} we give an example for $\gr(3,6)$. On this seed the cluster variables can be identified with Pl\"ucker coordinates $\pl{I}$ living in the coordinate ring of the affine cone. See \cite{scott-grassmannians,gekhtman2010cluster,cdfl-TableuxIndexing,fwz-IntroChap6} for full details.
\end{example}

\begin{definition}
    Given a mutable index $k$, we produce a new matrix $\mu_k(\exchangeMatFull{B})$ via \keyword{matrix mutation}:
    \begin{equation}
        \mu_k(\exchangeMatFull{B})_{i,j} = \begin{cases}
            -\exchangeMatFull{B}_{i,j} & i = k \text{ or } j = k\\
            \exchangeMatFull{B}_{i,j}+ \plus{\exchangeMatFull{B}_{i,k}}\plus{\exchangeMatFull{B}_{k,j}} - \minus{\exchangeMatFull{B}_{i,k}}\minus{\exchangeMatFull{B}_{k,j}} & \text{otherwise}
        \end{cases}
    \end{equation}
    where $\plus{x} = \max(x,0)$ and $\minus{x} = \min(x,0)$.
\end{definition}

\begin{definition}
     A \keyword{(cluster) seed pattern} is an $n$ regular tree whose vertices are labeled by seeds and edges are labeled $1,\dots,n$. The seeds $(\exchangeMatFull{B},\vec{\xvar{x}}) \xleftrightarrow{k} (\exchangeMatFull{B}',\vec{\xvar{x}}')$ satisfy the relations that $\exchangeMatFull{B}' = \mu_k(\exchangeMatFull{B})$ and 
     \begin{equation}
         \xvar{x}_{\ell}' = \mu_k(\xvar{x}_{\ell}) = \begin{cases}
             \xvar{x}_\ell & \ell \neq k \\
             \frac{1}{\xvar{x}_k}\left(\prod\limits_{i \st B_{ki}<0} \xvar{x}_i^{|B_{ki}|} + \prod\limits_{j \st B_{kj}>0}\xvar{x}_j^{B_{kj}} \right) & \ell = k
         \end{cases}
     \end{equation}
     Given a seed $S_0 = (\exchangeMatFull{B},\vec{\xvar{x}})$ we produce a seed pattern by starting at $S_0$ and performing all possible sequences of mutations. Note that $\mu_k(\mu_k(S)) = S$ so it makes sense to consider the edges of the $n$ regular tree to be unoriented.
\end{definition}
\begin{definition}
    A \keyword{$Y-$seed pattern} is similarly an $n$ regular tree with vertices labeled by $Y-$seeds. The key difference is the relation between two seeds on an edge, $(\exchangeMatFull{B},\vec{\yvar{y}}) \xleftrightarrow{k} (\exchangeMatFull{B}',\vec{\yvar{y}}')$. The exchange matrices are sill related by matrix mutation, but the variables are now related by \keyword{$Y$-mutation}:
    \begin{equation}
        \yvar{y}'_i = \mu_k(\yvar{y}_i) = \begin{cases}
            \yvar{y}_k^{-1} & i = k \\
            \yvar{y}_i(1+\yvar{y}_k^{-1})^{-b_{ik}} & b_{ik}>0\\
            \yvar{y}_i(1+\yvar{y}_k)^{-b_{ik}} & b_{ik} \leq 0\\
        \end{cases}
    \end{equation}
\end{definition}
\begin{remark}
    This definition of $Y-$variable mutation agrees with \cite{ahsl-ClusterConfigurationOfFiniteType} and \cite{fg-ClusterEnsembles}. However it is the opposite convention to \cite{fwz-IntroChap1-3} whose $Y-$variables correspond to the inverse of the $Y-$variables defined above.
\end{remark}

\begin{definition}
    The \keyword{cluster algebra} $\clusterAlg{A}_{Q}$ generated by a seed $(Q,\vec{\xvar{x}})$ is the subalgebra of the ring of rational functions on $\vec{\xvar{x}}$ generated by all cluster variables obtained by performing all possible mutations. Let $\xindexSet$ be the set indexing all mutable cluster variables and $\xindexSet_f$ be the set of frozen variables. Then $\clusterAlg{A} = \Z[x_i \st i \in \xindexSet_f][\xgamma\st \root{\gamma}\in \xindexSet]/I$ where $I$ is the ideal generated by the $X-$exchange relation for every possible mutation.\\ 
    One similarly defines the \keyword{$Y-$algebra}, $\yAlg{A}_Q = \Z[\yvar{y}]/J$ where $\yvar{y}$ ranges over every $Y-$variable in the $Y-$seed pattern and $J$ is generated by $Y-$exchange relations.
\end{definition}
\begin{remark}
    Matrix mutation preserves the rank of $\exchangeMatExtended{B}$ and thus we consider full rank a property of a cluster algebra/seed pattern.  
\end{remark}

\begin{claim}
    Let $\A$ and $\yAlg{Y}$ be cluster and $Y-$algebras for the same initial seed. There is a map $p\colon \yAlg{Y} \rightarrow \A$ given on each seed by 
    \begin{equation}
        \yvar{y}_i \mapsto \prod x_j^{\exchangeMatExtended{B}_{ij}}.
    \end{equation}
\end{claim}
\begin{proof}
    To check the map is well defined one must verify it commutes with $X$ and $Y$ mutation. See  \cite{gsv-ClusterAlgebraPoissonGeometry},\cite{fg-ClusterEnsembles} or \cite{fwz-IntroChap1-3} for the details of the proof.
\end{proof}

\begin{definition}
    The \keyword{localized cluster algebra} $\localizedClusterAlg{A} $ is obtained by localizing the set of mutable cluster variables.
\end{definition}
\begin{definition}
    The \keyword{$\semifield{R}$-points} of a cluster algebra in a semifield $\semifield{R}$ is the set $\Hom(\localizedClusterAlg{A},R)$. Each element can be identified with a choice of element in $\semifield{R}$ for each cluster variable satisfying the exchange relations.
\end{definition}

\begin{example}
    In type $A_1$ there are only two distinct seeds. There are two mutable cluster variables $\xvar{x}_{-\root{\gamma}}$ and $\xvar{x}_{\root{\gamma}}$. Thus $\A = \Z[\xvar{x}_{-\root{\gamma}},\xvar{x}_{\root{\gamma}}]/\ideal{\xvar{x}_{-\root{\gamma}}\xvar{x}_{\root{\gamma}} = 2}$. In \Cref{fig:A1RpointNoFrozen} we plot the real points. The positive points are the hyperbola in the first quadrant. In this somewhat trivial case the algebra is already localized at $x_{\root{\gamma}}$ and $x_{-\root{\gamma}}$ so the real points of $\A$ and $\localizedClusterAlg{A}$ are the same.\\
    In \Cref{fig:A1RpointsFullRank} we see a contour plot for the $A_1$ cluster algebra with one frozen variable $f_1$. We note that the frozen variable allows $x_{-\root{\gamma}}$ and $x_{\root{\gamma}}$ to take values at any pair of positive real numbers.

    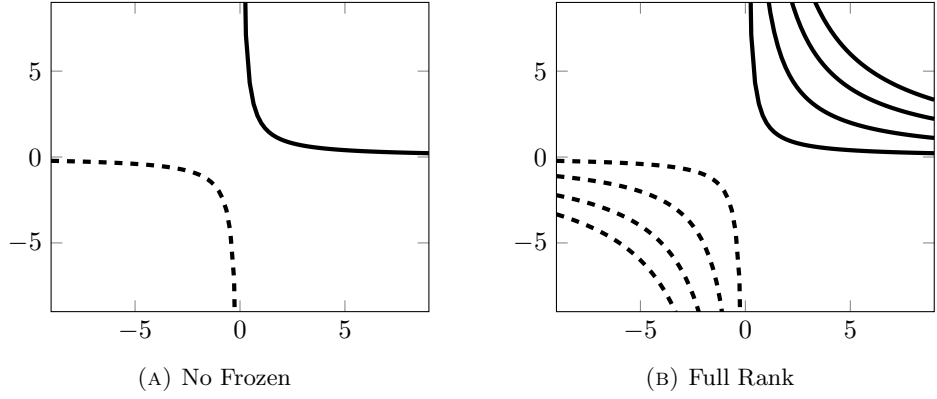
\begin{figure}
        \centering
        \begin{subfigure}[]{0.4\textwidth}
        \centering
            \begin{tikzpicture}
            \begin{axis}[xmin=-9,xmax=9,ymin=-9, ymax=9, samples=50, width=\textwidth]
              \addplot[ultra thick,domain=0.1:9] (x,2/x);
              \addplot[dashed,  ultra thick,domain=-9:-0.1] (x,2/x);
            \end{axis}
        \end{tikzpicture}
        \caption{No Frozen}
        \label{fig:A1RpointNoFrozen}
        \end{subfigure}       
    \begin{subfigure}[]{0.4\textwidth}
        \centering
        \begin{tikzpicture}
            \begin{axis}[xmin=-9,xmax=9,ymin=-9, ymax=9, samples=50, width=\textwidth]
              \addplot[ultra thick,domain=0.1:9] (x,2/x);
              \addplot[dashed,  ultra thick,domain=-9:-0.1] (x,2/x);
              \addplot[ultra thick,domain=0.1:9] (x,10/x);
              \addplot[dashed,  ultra thick,domain=-9:-0.1] (x,10/x);
              \addplot[ultra thick,domain=0.1:9] (x,20/x);
              \addplot[dashed,  ultra thick,domain=-9:-0.1] (x,20/x);
              \addplot[ultra thick,domain=0.1:9] (x,30/x);
              \addplot[dashed,  ultra thick,domain=-9:-0.1] (x,30/x);
            \end{axis}  
        \end{tikzpicture}
        \caption{Full Rank}
        \label{fig:A1RpointsFullRank}
    \end{subfigure}

        \caption{Real Points of $A_1$ Cluster Algebra}
        \label{fig:A1Rpoints}
    \end{figure}
\end{example}

\begin{claim}
    If $R$ has the additional property that the sum of two nonzero elements is nonzero then for every seed $(Q,\vec{\xvar{x}})$ of $\clusterAlg{A}$ every map $\Z[\vec{\xvar{x}}] \rightarrow R$ extends to a unique map $\localizedClusterAlg{A} \rightarrow{R}$.\\
    In particular, one can specify an $\R_{>0}$-point of a cluster algebra by choosing a positive number for each cluster variable in any seed.
\end{claim}
We call the $\R_{>0}$ points of a cluster algebra the \keyword{positive points} and write $\pospoints(\A)$ or $\pospoints$ if the choice of $\A$ is understood from the context. 

\subsection{Finite Type}
Let $D$ be a Dynkin diagram of finite type. 
\begin{definition}
    A quiver is of \keyword{Dynkin type} if its mutable part is an orientation of a Dynkin diagram with every node either a source or a sink. We say a cluster algebra $\A$ is of type $D$ if $\A$ contains a seed with quiver isomorphic to $D$. A \keyword{Dynkin seed} is a seed of $\A$ whose quiver is of Dynkin type. 
\end{definition}
\begin{proposition}
    Let $D$ be a connected Dynkin diagram and $h$ the associated Coxeter number. The cluster algebra of type $D$ contains $h+2$ Dynkin seeds. Each mutable cluster variable appears at a unique source and a unique sink in this set of seeds. 
\end{proposition}
\begin{proof}
    This follows from \cite{fz-ClusterII-finiteType}. As a consequence the set of mutable cluster variables is in bijection with set of almost positive roots $\Phi_{\geq -1}$ of the associated root system. 
\end{proof}

\begin{remark}
    The set of Dynkin seeds can be obtained by starting at any Dynkin seed and then repeatedly mutating at the full set of sources\cite{fz-ClusterII-finiteType}. This sequence seeds can be assembled into the Auslander-Reiten quiver of the associated cluster category (\cite{keller-clusterCategories}).
\end{remark}

While there are typically many more $Y-$variables than cluster variables, in a finite cluster algebra there is a canonical choice of $Y-$variable for each cluster variable, $\yvar{y}_\root{\gamma}$. We take $\yvar{y}_\root{\gamma}$ to be the $Y-$variable at the node associated to $\xvar{x}_\root{\gamma}$ in the Dynkin seed where $\xvar{x}_\root{\gamma}$ is a source. These $Y-$variables were used by Fomin and Zelivnsky to prove Zamolodchikov's conjecture \cite{fz-Ysystems}. As a consequence of this proof they established a Laurent phenomena analogous to the Laurent phenomena for cluster variables for this subset of $Y-$variables. In particular
\begin{theorem}\label{thm:YLaurent}
    Every $\yvar{y}_{\root{\gamma}}$ is of the form $\frac{N_\root{\gamma}}{\vec{p}^\root{\gamma}}$ where $N_\root{\gamma}$ is a polynomial in $p_i$ with constant term 1 and $\vec{p}^\vec{\root{\gamma}} = \prod p_i^{\root{\gamma}_i}$ and $\root{\gamma} = \sum \gamma_i e_i$ for $\{e_i\}$ the set of simple roots.
\end{theorem}
\begin{proof}
    This is a rephrasing of Theorem 1.5 of \cite{fz-Ysystems}.
\end{proof}
\begin{remark}
    The parameters $p_i$ are exactly the value of the $Y-$ variables at negative simple roots $-e_i$. These roots occur in two distinct Dynkin seeds. Usually the algebra is parameterized by choosing initial values from a single cluster. To obtain the same formulas as above from a single Dynkin seed we take the sources to be $p_i$ and the sinks to be $p_j^{-1}$. Then mutation at each sink produces the $Y-$ variable $p_i$ associated to the root $-e_i$ as needed.
\end{remark}

\begin{example}
    In type $A_3$ there are 6 Dynkin seeds as seen in \Cref{fig:A3BipartiteBelt} containing 9 distinct cluster variables. The values of each cluster variable and $Y-$variable are given in \Cref{fig:A3variables}. The exchange matrix for the initial seed without any frozen variables is $\begin{bmatrix}
        0 & 1 &0\\
        -1 & 0 & -1 \\
        0 & 1 & 0
    \end{bmatrix}$. This matrix is not full rank, it has null vector $(1,0,-1)$. However one can add frozen nodes to make the extended matrix full rank. For example an initial Dynkin seed for the cluster algebra associated to $\gr(2,6)$ is 
    \[ \begin{bmatrix}
        0& 1& 0& -1& 1& -1& 0& 0& 0\\
        -1& 0& -1& 1& 0& 0& 1& 0& 0\\
        0& 1& 0& 0& 0& 0& -1& 1& -1
    \end{bmatrix}.\]
    We observe the map from $Y-$variables to $X-$variables fails to be injective in the first case. Here the two $Y-$variables at the sources are: \[\yvar{y}_{[-1,0,0]} = \yvar{y}_{[0,0,-1]} = \xvar{x}_{[0,-1,0]}.\] In the Grassmannian case these variables are distinct: \[\yvar{y}_{[-1,0,0]} = \yvar{y}_{\pl{24}} = \frac{\pl{14}\pl{23}}{\pl{12}\pl{34}} \hspace{2pc}\text{and}\hspace{2pc}\yvar{y}_{[0,0,-1]} = \yvar{y}_{\pl{15}} = \frac{\pl{14}\pl{56}}{\pl{16}\pl{45}}.\]
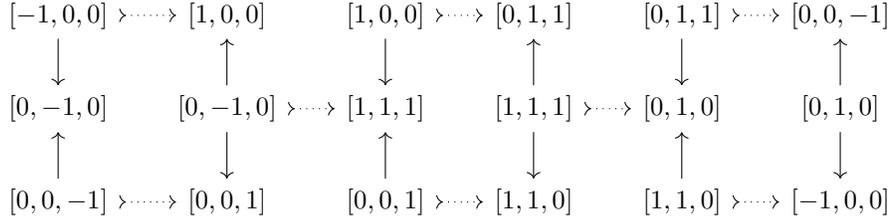
\begin{figure}[hb]
    \centering
    \begin{tikzcd}[column sep=1.5pc]
{[-1,0,0]} \arrow[d] \arrow[r, dotted, tail] & {[1,0,0]}                                              & {[1,0,0]} \arrow[d] \arrow[r, dotted, tail] & {[0,1,1]}                                             & {[0,1,1]} \arrow[d] \arrow[r, dotted, tail] & {[0,0,-1]}                                            &    \\
{[0,-1,0]}                                   & {[0,-1,0]} \arrow[d] \arrow[u] \arrow[r, dotted, tail] & {[1,1,1]}                                   & {[1,1,1]} \arrow[u] \arrow[d] \arrow[r, dotted, tail] & {[0,1,0]}                                   & {[0,1,0]} \arrow[d] \arrow[u]  \\
{[0,0,-1]} \arrow[u] \arrow[r, dotted, tail] & {[0,0,1]}                                              & {[0,0,1]} \arrow[u] \arrow[r, dotted, tail] & {[1,1,0]}                                             & {[1,1,0]} \arrow[u] \arrow[r, dotted, tail] & {[-1,0,0]}                                            &   
\end{tikzcd}

    \caption{Auslander-Reiten walk in type $A_3$}
    \label{fig:A3BipartiteBelt}
\end{figure}

\begin{figure}[!hb]
    \centering
${\renewcommand{\arraystretch}{1.5}
\begin{array}{r||c c | c |  }
     \text{Roots} & [-1,0,0] &  [0,0,-1] & [0,-1,0] \\
     \text{Cluster Variables}& x_1 & x_3 & x_2 \\
     \text{$Y-$Variables} & p_1 & p_3 & p_2 \\
\hline\hline
     \text{Roots}  & [1,0,0] & [0,0,1] & [1,1,1]  \\
     \text{Cluster Variables}& \frac{1+x_2}{x_1} & \frac{1+x_2}{x_3} & \frac{1+x_1x_3+2x_2+x_2^2}{x_1x_2x_3}\\
     \text{$Y-$Variables} & \frac{1+p_2}{p_1} & \frac{1+p_2}{p_3} & \frac{1+p_1+2p_2+p_3+p_1p_2+p_1p_3+p_2p_3+p_2^2}{p_1p_2p_3}\\
\hline \hline
     \text{Roots} &  [0,1,1] & [1,1,0] & [0,1,0] \\
     \text{Cluster Variables}& \frac{1+x_1x_3+x_2}{x_2x_3} & \frac{1+x_1x_3+x_2}{x_1x_2} & \frac{1+x_1x_3}{x_2}\\
     \text{$Y-$Variables} & \frac{1+p_1+p_2+p_3+p_1p_3}{p_2p_3} & \frac{1+p_1+p_2+p_3+p_1p_3}{p_1p_2} & \frac{1+p_1+p_3+p_1p_3}{p_2}
    \end{array}}$
    \caption{Cluster Variables and $Y-$Variables in $A_3$}
    \label{fig:A3variables}
\end{figure}

\end{example}

\begin{definition}
    The \keyword{compatibility degree} of two roots $\root{\gamma},\root{\omega}$ written $\compDegree{\root{\gamma}}{\root{\omega}}$ is the coefficient of $\root{\gamma}$ when $\root{\omega}$ is expressed as a linear combination of simple roots that includes $-\root{\gamma}$.\\
    The \keyword{compatibility degree} of two cluster variables $\compDegree{\xvar{x}_\root{\omega}}{\xvar{x}_\root{\gamma}}$ is the compatibility of the two roots. When $\root{\gamma} \neq \root{\omega}$, this corresponds to the power of $\xvar{x}_\root{\omega}$ in the denominator of the Laurent polynomial expression of $\xvar{x}_\root{\gamma}$ using the Dynkin seed containing $\xvar{x}_\root{\omega}$ as initial seed\cite{cp-denominatorCompatilityDegree}.
\end{definition}

\begin{example}
    The cluster algebra of type $C_2$ is full rank with no frozen variables. It has six Dynkin seeds containing 6 distinct cluster variables. If we take the initial seed to be $\xvar{x}_1 \rightarrow \xvar{x}_2$ with $\xvar{x}_1$ on the large node the variables are:
    \[\begin{aligned}
         \xvar{x}_1 = \xvar{x}_{[-1,0]} = \xvar{x}_1 \hspace*{2pc}& \xvar{x}_3=\xvar{x}_{[1,0]}=\frac{1+\xvar{x}_2}{\xvar{x}_1} \hspace*{0pc}& \xvar{x}_5=\xvar{x}_{[1,1]} = \frac{1+\xvar{x}_2+\xvar{x}_1^2}{\xvar{x}_1\xvar{x}_2} \\
         \xvar{x}_2=\xvar{x}_{[0,-1]} = \xvar{x}_2 \hspace*{2pc}& \xvar{x}_4=\xvar{x}_{[2,1]} = \frac{(1+\xvar{x}_2)^2+\xvar{x}_1}{\xvar{x}_1^2\xvar{x}_2}\hspace*{0pc}& \xvar{x}_6=\xvar{x}_{[0,1]} = \frac{1+\xvar{x}_1^2}{\xvar{x}_2}
    \end{aligned}
    \]
    The variables with odd indices are associated to the node with weight 2, while the variables with weight 1 are associated to the node with weight 1. 
    We compute the compatibility degree of $\xvar{x}_1 = \xvar{[-1,0]}$ with the other five roots:
    \begin{equation*}
        \compDegree{\xvar{x}_1}{\xvar{x}_2} = 0 \hspace{1pc} \compDegree{\xvar{x}_1}{\xvar{x}_3} = 1 \hspace{1pc} \compDegree{\xvar{x}_1}{\xvar{x}_4} = 2 \hspace{1pc} \compDegree{\xvar{x}_1}{\xvar{x}_5} = 1 \hspace{1pc}\compDegree{\xvar{x}_1}{\xvar{x}_6} = 0
    \end{equation*}
    Note that the compatibility degree is only symmetric for simply-laced Dynkin diagrams. This is manifest when comparing cluster variables associated to nodes with different weights. For example $\compDegree{\xvar{x}_4}{\xvar{x}_1} = 1 \neq 2= \compDegree{\xvar{x}_1}{\xvar{x}_4}$.
\end{example}

\section{Bounded Ratios}\label{sec:BoundedRatios}
\begin{definition}
    The \keyword{ratio space} $V^r \subset \localizedClusterAlg{A}$ is the set of Laurent monomials in mutable and frozen variables. We identify this space with an integer lattice of a real vector space $V$ of dimension $N+m$ where $N$ is the number of mutable cluster variables in the cluster algebra. The integer vector $(v_1,\dots,v_{N+m})$ corresponds to the ratio $\vec{\xvar{x}}^{\vec{v}} = \prod \xvar{x}_i^{v_i}$.
\end{definition}

We are interested in the subset of bounded ratios. 
\begin{definition}
    The \keyword{cone of bounded ratios}, $\cone{C}$, is subset of $V$ where the corresponding ratio is bounded over positive points of the cluster algebra. Formally \[\cone{C} = \spn_{\mathbb R^+}\{\ratiovec{v}\in V^r \st \exists L<\infty:  \forall p\in \pospoints : \vec{\xvar{x}}^\ratiovec{v}(p) < L \}\]
\end{definition}

\begin{example}\label{ex:A1ratios}
    Consider the cluster algebra for $A_1$ with one frozen variable (\Cref{fig:A1seedFullRank}). We identify $V$ with $\R^3$ using the basis $\{\xvar{x}_{-\root{\gamma}}, \xvar{x}_{\root{\gamma}}, f_1\}$. We translate several vectors into their corresponding ratios.
    \begin{align*}
        \begin{aligned}
            (-1,-1,0) =~& \xvar{x}_{-\root{\gamma}}^{-1}\xvar{x}_{\root{\gamma}}^{-1} = \frac{1}{1+\xvar{f}_1}\\
            (1,0,1) =~& \xvar{x}_{-\root{\gamma}}\xvar{f}_1
        \end{aligned}
        \hspace{2pc}
        \begin{aligned}
            (-1,-1,1)=~&\xvar{x}_{-\root{\gamma}}^{-1}\xvar{x}_{\root{\gamma}}^{-1}\xvar{f}_1^{1} = \frac{\xvar{f}_1}{1+\xvar{f}_1}\\
            (1,1,0) =~& \xvar{x}_{-\root{\gamma}}^1\xvar{x}_{\root{\gamma}}^1= 1+\xvar{f}_1
        \end{aligned}
    \end{align*}
    We observe the vectors in the first row are bounded and thus in $C$. The vectors in the second row are unbounded, and thus not in $C$.
\end{example}

\subsection{Torus Action}
We now explore a simple necessary condition for a ratio to be bounded. 

Each vector $\vec{\alpha} \in \ker(\exchangeMatExtended{B})$ and $z\in \multgroup{R}$ induces an automorphism $\phi_{\vec{\alpha},z}$ of $\clusterAlg{\A}\tensor R$ that sends $\xvar{x}_i \mapsto z^{\alpha_i}\xvar{x}_i$ for each initial mutable and frozen cluster variable. By Lemma 5.3 in \cite{gekhtman2010cluster} when $\vec{\alpha} \in \ker(\exchangeMatExtended{B})$, the induced action on all other cluster variables is multiplication by a power of $z$ as well.  We write $\weight_\vec{\alpha}(x_\root{\gamma})$ for this power. The weight extends to the full ratio space in the natural way, $\weight_\vec{\alpha}(\ratiovec{v}) = \sum \weight_\vec{\alpha}(x_\root{\gamma}) v_\root{\gamma}$.\\
As $\phi_{\vec{\alpha},z}$ is an automorphism of $\clusterAlg{\A}\tensor R$, it induces an action on the $R$ points of the algebra. In particular, if $z \in \R_{>0}$, this action preserves the set of positive points.\\
Finally we compute the action of $\phi_{\vec{\alpha},z}$ on the pairing between any ratio vector $\vec{v}$ and $R-$point $p$:
\begin{equation}\label{eqn:TorusAction}
    \ratiovec{v}(\phi_{\vec{\alpha},z} p) = (\phi_{\vec{\alpha},z} \ratiovec{v})(p) =  z^{\weight_\vec{\alpha}(\ratiovec{v})}\ratiovec{v}(p)
\end{equation}
\begin{remark}
    This action can be viewed as an algebraic torus $T$ acting by characters on the ratio space. See \cite{ahsl-ClusterConfigurationOfFiniteType,ls-cohomologyClusterVarieties}.
\end{remark}

\begin{lemma}
    If there is a vector $\vec{\alpha}\in \ker(\exchangeMatExtended{B})$ such that $\weight_\vec{\alpha}(\ratiovec{v}) \neq 0$ then $\ratiovec{v}$ is not a bounded ratio.
\end{lemma}
\begin{proof}
    If $\weight_\vec{\alpha}(\ratiovec{v}) >0 $ then take $z=2$ otherwise take $z=\frac{1}{2}$. Then \[\ratiovec{v}(\phi_{\vec{\alpha},z}^t(p)) = (z^{\weight_{\vec{\alpha}}(\ratiovec{v})})^t\cdot \ratiovec{v}(p).\] As $t \rightarrow \infty$, $(z^{\weight_{\vec{\alpha}}(\ratiovec{v})})^t \rightarrow \infty$ and thus $\ratiovec{v}$ is unbounded.
\end{proof}

\begin{corollary}\label{thm:BoundedImpliesWeight0}
    A bounded ratio has weight $0$ for any vector $\vec{\alpha} \in \ker(\exchangeMatExtended{B})$.
\end{corollary}

\begin{example}\label{ex:A1torusAction}
In $A_1$ with one frozen variable the vector $\vec{\alpha} = (-1,0)$ generates the kernel of the exchange matrix (\Cref{fig:A1seedFullRank}). The induced map on the cluster algebra is
\begin{equation*}
    (\xvar{x}_{-\gamma}, \xvar{x}_{\gamma},\xvar{f}_1) \xmapsto{\phi_{\vec{\alpha},z}} (z^{-1}\xvar{x}_{-\gamma}, z^{1} \xvar{x}_{\gamma},z^{0} \xvar{f}_1)
\end{equation*}
and the induced weight vector is $(-1,1,0)$. We observe that the orbits of this action are the hyperbolas in \Cref{fig:A1RpointsFullRank}.\\
Furthermore from \Cref{thm:BoundedImpliesWeight0}, the cone of bounded ratios is contained in the subspace with equal powers of $\xvar{x}_{-\root{\gamma}}$ and $\xvar{x}_{\root{\gamma}}$. To illustrate this containment, we compute the weights of the ratios in \Cref{ex:A1ratios}:
\[\weight_{\vec{\alpha}}(-1,-1,0) = 0 \hspace{1pc} \weight_{\vec{\alpha}}(-1,-1,1) = 0 \hspace{1pc} \weight_{\vec{\alpha}}(1,0,1) = -1 \hspace{1pc} \weight_{\vec{\alpha}}(1,1,0) = 0 .\]
We observe the first two bounded ratios have zero weight as expected and the third unbounded ratio has nonzero weight. However the final ratio is unbounded despite having zero weight.
\end{example}

\begin{remark}
    In the Grassmannian cluster algebra \Cref{thm:BoundedImpliesWeight0} can be rephrased using the natural algebraic torus action of $(\multgroup{\C})^n$ on $\gr(k,n)$ \cite{boocher2008generators}. This \keyword{ST0 condition} states that each column index appears an equal number of times in the numerator and denominator. This corresponds to a character given by assigning a variable weight $\alpha_i$ where $\alpha_i$ is the number of times column $i$ appears when expressed as a polynomial in Pl\"ucker coordinates.
\end{remark}
Note that all exchange relations corresponding to standard initial cluster of $\gr(k,n)$ are either short Pl\"ucker relations 
\[ \pl{ikS}\pl{j\ell S} = \pl{ij S}\pl{k\ell S} + \pl{i\ell S}\pl{jk S}\]
 where $S$ is a $k-2$ element subset of $[1,n]$ and $i,j,k,\ell \notin S$, or are of the form
 $$\pl{ijkS}~\mu(\pl{ijkS})=\pl{ikfS}~\pl{ijdS}~\pl{jkeS}+\pl{ikdS}~\pl{ijeS}~\pl{jkfS},$$
 where $S$ is $k-3$ element subset of $[1,n]$, such that $i,d,j,e,k,f \notin S$. In each such relation both monomials in RHS have same multiplicities for each index i. Thus the vector $\alpha_i$ described above is in the kernel of the exchange matrix. 

\begin{lemma}\label{thm:LinearDimensionOfCone}
    When the extended exchange matrix has full rank, the cone of bounded ratios is contained in a space of dimension equal to number of cluster variables, $N$.
\end{lemma}
\begin{proof}
    Since the extended exchange matrix has full rank $n$,  $\dim(\ker(\exchangeMatExtended{B})) =m$. Thus for a ratio to have weight 0 with respect to the entire kernel, $\ratiovec{v}$ must satisfy $m$ linearly independent conditions. Thus the cone of bounded ratios is contained in a space of dimension $N+m -m = N$. 
\end{proof}

\subsection{U Variables}

In \cite{ahsl-ClusterConfigurationOfFiniteType} they define the cluster configuration space associated to a finite Dynkin diagram $D$. These spaces recover the classical configuration space of $n$ distinct points on the projective line, $M_{0,n}$, by taking $D$ to be a Dynkin diagram of type $A_n$. We recall the definitions here.
\begin{definition}
    The \keyword{cluster configuration algebra}, $\uAlg{U}_{D}$, associated to a Dynkin diagram $D$ is the ring $\Z[u_\root{\gamma}^{\pm 1}]/I_D$ where $I_D$ is the ideal generated by the equations of the form
    \begin{equation}\label{eqn:Usystem}
        \uvar{u}_\root{\gamma} + \prod_{\root{\omega}\neq \root{\gamma}} \uvar{u}_\omega^{\compDegree{\root{\omega}}{\root{\gamma}}} =1
    \end{equation}
    for each mutable cluster variable $\root{\gamma} \in \xindexSet$ of the associated cluster algebra $\clusterAlg{A}_D$.
\end{definition}

As with the cluster algebra, we can discuss $R$ points of the cluster configuration algebra. Following \cite{ahsl-ClusterConfigurationOfFiniteType} we write $\uPoints$ for the $\C$ points of $\uAlg{U}_D$ and $\uPosPoints$ for the $\R_{>0}$ points.

Surprisingly, the space $\uPosPoints$ is $n$ dimensional where $n$ is the number of nodes in $D$. In fact there is a map from $\pospoints$ to $\uPosPoints$ given by \begin{equation}\label{eqn:Uparam}
    (\xgamma)_{\gamma \in \xindexSet} \mapsto (v_{\gamma})_{\gamma \in \xindexSet} = \left(\frac{\yvar{y}_\gamma}{1+\yvar{y}_\gamma}\right) = \left(\frac{\prod \xomega}{\xgamma \xgamma'} \right) 
\end{equation}
where $\yvar{y}_{\root{\gamma}}$ is the $Y-$coordinate in the Dynkin seed at $\xgamma$ where $\root{\gamma}$ is a source. The final equality uses the map from $Y-$coordinates to $X-$coordinates and the product in the numerator is taken over all $\root{\omega}$ that are out neighbors of $\root{\gamma}$ and $\xgamma'$ is the $X-$coordinate obtained by mutating at $x_\gamma$ in this seed.
\begin{proposition}
    When $\clusterAlg{\A}$ is of full rank, the map given in \Cref{eqn:Uparam} is a bijection $\pospoints\sslash T \rightarrow \uPosPoints$. In fact this induces an isomorphism $\uAlg{U}_D \rightarrow \clusterAlg{A}_D^{T}$ given by $\uvar{u}_\root{\gamma} \mapsto v_\root{\gamma}$.
\end{proposition}
\begin{proof}
    This follows from Theorem 4.2 of \cite{ahsl-ClusterConfigurationOfFiniteType}. They describe this isomorphism on the level of rings and on the level of $\C$ points. However as the torus action restricts to the set of positive points, the isomorphism works in this context as well.
\end{proof}

\begin{example}
    The cluster configuration algebra of type $A_1$ is defined by the equations
    \begin{equation*}
        \uvar{u}_{-\root{\gamma}} + \uvar{u}_{\root{\gamma}} = 1 \hspace{2pc} \uvar{u}_{\root{\gamma}} + \uvar{u}_{-\root{\gamma}} = 1
    \end{equation*}
    In \Cref{fig:A1UPositivePoints} we see $\uPosPoints$. Here we can observe the requirement that $\A$ is of full rank. If we take the $A_1$ cluster algebra with no frozen nodes the map from $\pospoints$ to $\uPosPoints$ is of the form 
    \begin{equation*}
        (\xvar{x}_{-\gamma},\xvar{x}_{\gamma}) \mapsto \left(\frac{1}{\xvar{x}_{-\gamma}\xvar{x}_{\gamma}}, \frac{1}{\xvar{x}_{-\gamma}\xvar{x}_{\gamma}}\right) = \left(\frac{1}{2},\frac{1}{2}\right)
    \end{equation*}
    and thus doesn't fill the entire space of positive points. However when we add a frozen variable to make a full rank cluster algebra the map is now of the form:
    \begin{equation*}
        (\xvar{x}_{-\gamma},\xvar{x}_{\gamma},\xvar{f}_1) \mapsto \left(\frac{1}{\xvar{x}_{-\gamma}\xvar{x}_{\gamma}}, \frac{\xvar{f}_1}{\xvar{x}_{-\gamma}\xvar{x}_{\gamma}}\right) = \left(\frac{1}{1+\xvar{f}_1},\frac{\xvar{f}_1}{1+\xvar{f}_1}\right).
    \end{equation*}
    
    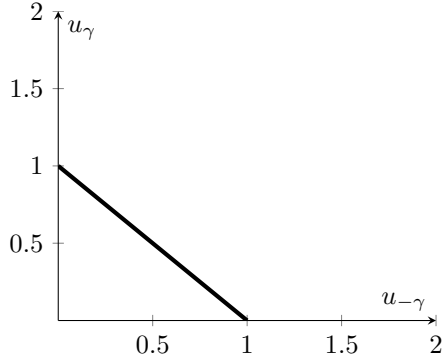
\begin{figure}[!hb]
        \centering
        \begin{tikzpicture}
            \begin{axis}[xmin=0,xmax=2,ymin=0, ymax=2, samples=50, width=.4\textwidth, axis x line = middle, axis y line = middle, xlabel={$u_{-\gamma}$}, ylabel={$u_{\gamma}$}]
                \addplot[ultra thick,domain=0:1] (x,1-x);
            \end{axis}
        \end{tikzpicture}
        \caption{Positive Points of the Cluster Configuration Algebra of type $A_1$}
        \label{fig:A1UPositivePoints}
    \end{figure}
    
\end{example}

\begin{example}\label{ex:C2uvars}
    The cluster configuration algebra of type $C_2$ is defined by the orbit of the equations
    \[ \uvar{u}_1 + \uvar{u}_3\uvar{u}_4\uvar{u}_5=1 \hspace{2pc} \uvar{u}_2+\uvar{u}_4\uvar{u}_5^2\uvar{u}_6 = 1\]
    under increasing the indices by $2$ modulo $6$. Under the map from \Cref{eqn:Uparam} the solutions are parameterized by
    \begin{align*}
        \uvar{v}_1 = \frac{\xvar{x}_2}{\xvar{x}_1\xvar{x}_3} \hspace{1pc} \uvar{v}_3 = \frac{\xvar{x}_4}{\xvar{x}_3\xvar{x}_5} \hspace{1pc} \uvar{v}_5 = \frac{\xvar{x}_6}{\xvar{x}_5\xvar{x}_1} \hspace{2pc}
        \uvar{v}_2 = \frac{\xvar{x}_3^2}{\xvar{x}_2\xvar{x}_4} \hspace{1pc} \uvar{v}_4 = \frac{\xvar{x}_5^2}{\xvar{x}_4\xvar{x}_6} \hspace{1pc} \uvar{v}_6 = \frac{\xvar{x}_1^2}{\xvar{x}_6\xvar{x}_2}.
    \end{align*}
\end{example}

\begin{claim}\label{thm:Ubounded}
    Each variable $u_\root{\gamma}$ is bounded over $\uPosPoints$.
\end{claim}
\begin{proof}
    This is clear from the form of \Cref{eqn:Usystem} over the positive reals, each $u_\root{\gamma}$ is $1$ minus a positive number and thus is bounded above by $1$.
\end{proof}

\begin{theorem}\label{thm:UboundaryLimits}
    When $\clusterAlg{A}$ is of full rank and finite type, there is a choice of initial cluster variables parameterized by $t\in \R_{>0}$ such that $\lim\limits_{t \rightarrow 0} \uvar{v}_{\root{\gamma}} = 0$  and for $\root{\omega}\ne\root{\gamma}$, $0< \lim\limits_{t\rightarrow 0}\uvar{v}_{\root{\omega}} <\infty$. 
\end{theorem}
\begin{proof}
    Consider the Dynkin seed where $\xgamma$ is a source at vertex $i$. We consider this seed to be the initial seed and thus can express all other cluster variables in terms of the initial variables $\xvar{x}_{1},\dots,\xvar{x}_{n},\xvar{x}_{n+1},\dots,\xvar{x}_{n+m}$.\\
    Assume the initial cluster variables are of the form $\xvar{x}_j = t^{\beta_j}$. Then the initial $Y-$variables are of the form $\yvar{y}_j =\prod \xvar{x}_k^{\exchangeMatExtended{B}_{jk}} = \prod t^{\exchangeMatExtended{B}_{jk}\beta_k}$. Thus the vector of powers of initial $Y-$variables as a function of $t$ is $\vec{\yvar{y}} = \exchangeMatExtended{B}\vec{\beta}$. As $\exchangeMatExtended{B}$ is of full rank, the equation $e_i = \exchangeMatExtended{B}\vec{\beta} $ has a solution. We use this $\beta$ to define the initial cluster variables. By construction, the other initial $Y-$variables are constant with respect to $t$. By \Cref{thm:YLaurent} all noninitial $Y-$variables are Laurent polynomials in variables $p_1,\dots, p_n$ whose numerator has constant term $1$. By our choice of initial cluster variables we have that $p_i = t$ and $p_j = 1$ for $j \neq i$, since $p_i = \yvar{y}_i$ and $p_j = \yvar{y}_j^{\pm 1}$ depending on if $j$ is a source or sink. We then have 
    \begin{equation}
        \lim_{t\rightarrow 0}\yvar{y}_\root{\omega} = \lim_{t\rightarrow 0} \frac{C_\root{\omega} +o(t)}{t^{\compDegree{\root{\gamma}}{\root{\omega}}}} = \begin{cases}
            C_\root{\omega} & \compDegree{\root{\gamma}}{\root{\omega}} = 0\\
            \infty & \compDegree{\root{\gamma}}{\root{\omega}} > 0
        \end{cases} 
    \end{equation}
    with $C_\root{\omega}$ a positive constant. As $\uvar{v}_\root{\omega} = \frac{\yvar{y}_\root{\omega}}{1+\yvar{y}_\root{\omega}}$ we see that $\lim\limits_{t \rightarrow \infty} \uvar{v}_\root{\omega}$ is a positive constant if $\root{\omega}$ and $\root{\gamma}$ are compatible and equal to $1$ otherwise. 
\end{proof}

\begin{corollary}\label{thm:Uindependence}
    The ratios $\uvar{v}_\root{\gamma}$ correspond to independent vectors in the root space. In other words there is no vector of powers $\vec{\lambda}$ such that  $\uvar{v}_\root{\gamma} = \prod\limits_{\root{\omega}\neq \root{\gamma}} \uvar{v}_\root{\omega}^{\lambda_\root{\omega}}$.
\end{corollary}
\begin{proof}
    Assume for contradiction that such a product existed. At least one power $\lambda_\omega$ is nonzero, otherwise $\uvar{v}_\root{\gamma}$ would be identically $1$ for all choices of initial cluster variables. However by \Cref{thm:UboundaryLimits} there is a choice of initial seed with $\uvar{v}_\root{\gamma}$ going to $0$ at $t$ goes to $0$. \\
    Now consider the choice of initial seed given by \Cref{thm:UboundaryLimits} such that only $\uvar{v}_\root{\omega}$ goes to 0 while each other variable is bounded away from $0$. Since the set of ratios is finite there is a global lower bound $\lim\limits_{t\rightarrow\infty} \uvar{v}_{\root{\eta}} > L>0$ on all other ratios. Note that as all ratios are also bounded above by 1 (\Cref{thm:Ubounded}) the other product terms stay finite. If $\lambda_\root{\omega} >0$ then the entire product is sent to 0. On the other hand $\uvar{v}_\root{\gamma}$ is bounded below by $L$, contradicting the equality. Similarly if $\lambda_\root{\omega} <0$ the product is sent to infinity contradicting the fact that  $v_{\root{\gamma}}$ is bounded (\Cref{thm:Ubounded}).
\end{proof}

\begin{theorem}\label{thm:GeneratorsOfFullConeOfBoundedRatios}
    The cone of bounded ratios $C$ is generated by the set $\uvar{v}_\root{\gamma}$.
\end{theorem}
\begin{proof}
    From \Cref{thm:LinearDimensionOfCone} we know that $C$ is contained in a vector space of dimension $N$, the space of weight 0 ratios. \Cref{thm:Uindependence} states the the set $\{\uvar{v}_\root{\gamma}\}$ is a set of $N$ linearly independent ratios. As each $\uvar{v}_\root{\gamma}$ has weight 0, this gives a basis of the subspace of weight 0 ratios. Therefore every ratio of weight 0 has the form $\sum \lambda_\root{\gamma} v_\root{\gamma}$ for some vector $\vec{\lambda}$. We will show that if $\vec{\lambda}$ contains any negative entries the corresponding ratio is unbounded. \\
    Assume $\vec{\lambda}$ is a vector with  $\lambda_\root{\gamma}<0$. Then consider the ray from \Cref{thm:UboundaryLimits} sending $v_\root{\gamma}$ to 0 while keeping all other $v_\root{\omega}>L$ for some $L>0$. Then $v_\root{\gamma}^{\lambda_\gamma}$ goes to infinity while the rest of the ratio is bounded away from 0. Thus the whole ratio is unbounded along this ray and thus doesn't belong to the cone of bounded ratios. Therefore the cone of bounded ratios is exactly the positive linear combinations of $v_\root{\gamma}$ as claimed.
\end{proof}

\begin{corollary}
    The cone of bounded ratios in the cluster configuration algebra for a finite Dynkin diagram $D$ is the positive orthant.
\end{corollary}
\begin{proof}
    We can always find a cluster algebra of type $D$ with full rank. In this case the $\uvar{v}_\root{\gamma}$ exactly parameterize the $\uvar{u}_{\root{\gamma}}$.
\end{proof}
\begin{corollary}
    A ratio in cluster variables is bounded if and only if it is bounded by 1. If such a ratio has integer coefficients when expressed in the basis $\uvar{v_\root{\gamma}}$ it is also  subtraction free i.e. denominator minus numerator is a polynomial with no negative signs in cluster variables (not necessarily from the same cluster). 
\end{corollary}
\begin{proof}
    The extreme rays $\uvar{v}_\root{\gamma}$ clearly have both properties from \Cref{eqn:Uparam}. By the cluster relation we see the denominator minus numerator is simply the product of frozen variables coming into node $\root{\gamma}$. This inductively implies every integer linear combination is subtraction free as well by the following computation.
    \[\frac{a}{b}\frac{c}{d} \mapsto bd-ac = b(d-c)+c(b-a)\]
\end{proof}
\begin{example}\label{ex:BoundedButNotSubtractionFree}
    There are bounded ratios of integer powers of cluster variables that are not subtraction free.  
    In \Cref{ex:C2uvars} we computed the extreme rays of the cone of bounded ratios in type $C_2$. The product $\sqrt{\uvar{v}_2\uvar{v}_4\uvar{v}_6} = \frac{\xvar{x}_1\xvar{x}_3\xvar{x}_5}{\xvar{x}_2\xvar{x}_4\xvar{x}_6}$ is an integer combination of cluster variables. However if we express this ratio in the initial cluster $\xvar{x}_1 \rightarrow \xvar{x}_2$ we obtain  
    \begin{align*}
    x_2x_4x_6-x_1x_3x_5 =~&\xvar{x}_2 \left( \frac{(1+\xvar{x}_2)^2+\xvar{x}_1}{\xvar{x}^2_1\xvar{x}_2}\right)\left( \frac{1+\xvar{x}^2_1}{\xvar{x}_2}\right) - \xvar{x}_1 \left(\frac{1+\xvar{x}_2}{\xvar{x}_1} \right)\left( \frac{1+\xvar{x}_2+\xvar{x}^2_1}{\xvar{x}_1\xvar{x}_2}\right) \\
    =~& \frac{(1 + 2 \xvar{x}_1^2 + \xvar{x}_1^4 + 2 \xvar{x}_2 + 2 \xvar{x}_1^2 \xvar{x}_2 + \xvar{x}_2^2 + \xvar{x}_1^2 \xvar{x}_2^2) - \left( \xvar{x}_1 + \xvar{x}_1^3 + 2 \xvar{x}_1 \xvar{x}_2 + \xvar{x}_1^3 \xvar{x}_2 + \xvar{x}_1 \xvar{x}_2^2\right)}{x_1^2x_2} 
    \end{align*}
    If this ratio were subtraction free in cluster variables, then specializing to the initial cluster should produce a Laurent polynomial with positive coefficients. Since this didn't happen, the ratio cannot be subtraction free in cluster variables as claimed.  
   
\end{example}
o\begin{lemma}\label{thm:SufficientSubtractionFree}
    Consider the $(N+m) \times N$ matrix $U$, where $U_{i,j}$ is the power of the cluster variable $\xvar{x}_i$ in the extreme ratio $\uvar{v}_{\root{\gamma}_j}$. If there is a subset of rows $I$ such that $\det(U_I) = \pm 1$ then every bounded ratio is uniquely expressed as a positive integer combination of the $\uvar{v}_{\root{\gamma}}$.
\end{lemma}
\begin{proof}
    Consider a bounded ratio represented by the vector $\vec{v}$. By construction, the solution to the equation $U\vec{\lambda} = \vec{v}$ provides the coefficients $\lambda$ needed to express $\vec{v}$ as a combination of $\uvar{v}_{\root{\gamma}}$. If we restrict to the $I$ rows we see $U_I \vec{\lambda} = v_I$. Furthermore by assumption $\det(U_I) = \pm 1$ and thus by Cramer's rule $U_I^{-1}$ is an integer matrix. Therefore $\vec{\lambda} = U_I^{-1}v_I$ is integer as needed. Since the set of $\vec{v}_{\root{\gamma}}$ form a basis, $\vec{\lambda}$ is unique. 
\end{proof}
\begin{corollary}\label{thm:Gr3678SubtractionFree}
    The bounded ratios on the cluster algebras for $\gr(3,6), \gr(3,7)$ and $\gr(3,8)$ are subtraction free.
\end{corollary}
\begin{proof}
    All three of these cluster algebras satisfy the condition of \Cref{thm:SufficientSubtractionFree} by explicit computation.
\end{proof}

\begin{example}
    We now finish our running example of the $A_1$ cluster algebra with one frozen variable. From \Cref{eqn:Uparam} we have the two generating ratios are
    \[ \uvar{v}_{-\root{\gamma}} = \frac{1}{\xvar{x}_{-\root{\gamma}}\xvar{x}_\root{\gamma}} \hspace{2pc} \uvar{v}_{\root{\gamma}}=\frac{f_1}{\xvar{x}_{-\root{\gamma}}\xvar{x}_\root{\gamma}}\]
    From \Cref{ex:A1torusAction} we know the cone of bounded ratios is contained in the subspace where the power of $\xvar{x}_{-\root{\gamma}}$ equals the power of $\xvar{x}_{\root{\gamma}}$. Thus we can draw the cone of bounded ratios in \Cref{fig:A1BoundedCone}.
    \begin{figure}[!hb]
        \centering
        \begin{tikzpicture}
            \begin{axis}[xmin=-9,xmax=9,ymin=-9, ymax=9, samples=50, width=.4\textwidth, axis x line = middle, axis y line = middle, xlabel={$x_\gamma$}, ylabel={$f_1$}]
                \addplot[name path=A,ultra thick,domain=0:9] (-x,0);
                \addplot[name path=B,ultra thick,domain=0:9] (-x,x);
                \addplot[color=black,fill=pink,fill opacity=0.4]fill between[of=A and B, soft clip={domain=-9:0}];
            \end{axis}
        \end{tikzpicture}
        \caption{Cone of Bounded Ratios in $A_1$}
        \label{fig:A1BoundedCone}
    \end{figure}
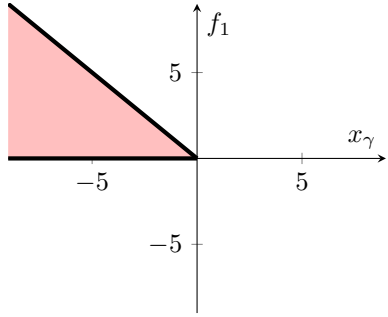
    Furthermore every bounded ratio is subtraction free by \Cref{thm:SufficientSubtractionFree} since the matrix $U$  is 
    \[ \begin{bmatrix} -1 & -1 \\ -1 & -1\\ 0 & 1\end{bmatrix}\].
    
\end{example}

\section{Bounded Ratios on Subsets}\label{sec:SubsetRatios}
The original study of totally nonnegative matrices was concerned with inequalities of determinants of matrix minors. When translated to Grassmannians, $\gr(k,n)$, these are inequalities among the Pl\"ucker coordinates. The most basic such ratio is the \emph{primitive ratio} \cite{soskin2023bounded,boocher2008generators} of the form 
\[ \frac{\pl{i(j+1)S}\pl{j(i+1) S}}{\pl{ijS}\pl{(i+1)(j+1)S}}\]
where $i<i+1<j<j+1$ are cyclically ordered and $S$ is a $k-2$ element subset of $[n]$ disjoint from $i,i+1,j,j+1$. It was conjectured in \cite{fallat2003multiplicative} that all bounded ratios on the Grassmannian can be factored into positive integer combinations of primitive ratios. However in $\gr(4,8)$, \cite{soskin2023bounded} found extreme ratios that are not primitive ratios. \\
Using the results of this paper we prove that the bounded ratios in Pl\"ucker coordinates on Grassmannians with finite type cluster structure ($\gr(2,n), \gr(3,6), \gr(3,7), \gr(3,8)$) are generated by primitive ratios. To accomplish this we give a concrete algorithm to reduce the full cluster bounded cone to a bounded cone on a subset of the cluster variables.

\begin{definition}
    Let $S \subset \xindexSet$. The \keyword{$S-$ratio space}, $V^{r}_S$, is the subset of $V^r$ generated by the variables in $S$. The \keyword{cone of bounded $S-$ratios}, $C_S$, is the cone of bounded ratios in $V_S$, where $V_S$ is a linear subspace of $V$ spanned by $V^{r}_S$.
\end{definition}
\begin{proposition}\label{thm:ExtremeRaysOfSBoundedCone}
    The cone of bounded $S-$ratios $C_S$ is the intersection of the full cone of bounded ratios with $V_S$. Moreover, when $\clusterAlg{A}$ is of full rank, the extreme rays of $C_S$ are the minimal positive linear combinations of the extreme rays of the full cone that are contained in $V_S$.
\end{proposition}
\begin{proof}
    Since every bounded $S-$ratio is a usual bounded ratio the first statement is clear. In the full rank case every ratio is uniquely expressed as product of $v_{\root{\gamma}}$. Thus if no strict subproduct lands in $V_S$ the corresponding ratio cannot be factored in $C_S$.
\end{proof}

Given \Cref{thm:ExtremeRaysOfSBoundedCone} it is then straightforward to compute the extreme rays of $C_S$. Let $U$ be the matrix whose columns are indexed by $\xindexSet$ and rows are indexed by the full set of cluster variables without $S$ with $U_{i,j}$ the power of $\xvar{x}_{i}$ in the extreme ratio $\uvar{v}_j$. 
\begin{claim}
    Let $\vec{\lambda} \in \ker(M)$. Then $\prod v_{\root{\gamma}}^{\lambda_{\root_\gamma}}$  is in  $V_S$. Furthermore if every element of $\vec{\lambda}$ is nonnegative the corresponding ratio is in $C_S$.
\end{claim}
\begin{proof}
    By construction the entries of $M\vec{\lambda}$ are the power of corresponding the cluster variable not in $S$ in the resulting ratio. Thus if $\vec{\lambda} \in \ker(M)$ all variables not in $S$ vanish. If the coefficients of $\vec{\lambda}$ are all nonnegative the ratio is bounded and thus lies in $C_S$.
\end{proof}
Therefore $C_S$ is defined by the equations $M\vec{\lambda} = 0$ and $\lambda_{\root{\gamma}} \geq 0$. Such a system of linear inequalities can be easily solved by programs like  Normaliz \cite{Normaliz}. 

We now return our focus to the Grassmannian and specialize $S$ to the set of Pl\"ucker variables. In $\gr(2,n)$ every cluster variable is a Pl\"ucker coordinate. Thus the full cone is the same as the Pl\"ucker cone and we have the following:
\begin{proposition}\label{thm:Gr2nFactor}
    A ratio of Pl\"ucker coordinates on $\gr(2,n)$ is bounded if and only if it can be factored into a positive integer combination of primitive ratios. 
\end{proposition}
\begin{proof}
    The cluster structure on $\gr(2,n)$ is of type $A_{n-3}$. We can see the generator $v_{ij}$ given by  \Cref{thm:GeneratorsOfFullConeOfBoundedRatios} is the primitive ratio indexed by $\{(i-1),(j-1)\}$ since the Dynkin seeds correspond to the zig-zag triangulations of an $n$-gon.  
\end{proof}

\begin{example}
    Consider the cluster algebra $\gr(2,n)$. Let $S$ be the set of Pl\"ucker coordinates with indices in $\{1,\dots,n-1\}$, i.e. the set of all cluster variables in $\gr(2,n-1)$. Using \Cref{thm:ExtremeRaysOfSBoundedCone} we can compute the extreme rays of $\gr(2,n-1)$ which we call $\uvar{u}_{ij}(n-1)$ using the extreme rays of $\gr(2,n)$, $\uvar{u}_{ij}(n)$. Explicitly we have
    \[ \uvar{u}_{1j}(n-1) = \uvar{u}_{1j}(n)\cdot \uvar{u}_{jn} \hspace{2pc} \uvar{u}_{ij}(n-1) = \uvar{u}_{ij}(n) \text{ otherwise } \]
    This follow from the following simple computation
    \[ \uvar{u}_{1j}(n-1) = \frac{\pl{j(n-1)}\pl{1(j-1)}}{\pl{1j}\pl{(j-1)(n-1)}} = \frac{\pl{1(j-1)}\pl{jn}}{\pl{1j}\pl{(j-1)n}}\frac{\pl{j(n-1)}\pl{(j-1)n}}{\pl{jn}\pl{(j-1)(n-1)}}\]
    and the observation that $\uvar{u}_{ij}(n) \in V_S$ if $1<i$ and $j < n$.
\end{example}

This inductive structure of the extreme rays allows us to strengthen \Cref{thm:Gr2nFactor} to the following theorem:
\begin{theorem}\label{thm:Gr2nSubtractionFree}
    Every bounded ratio of Pl\"ucker coordinates on $\gr(2,n)$ can be factored as an integer combination of primitive ratios. As such every such bounded ratio is subtraction free and bounded by 1.
\end{theorem}
\begin{proof}
    By \Cref{thm:SufficientSubtractionFree} it suffices to find a minor of the matrix of $U-$variables with determinant $\pm1$. In \cite[Theorem 4.2]{soskin2023bounded} it is proved that such a minor exists for any $\gr(k,n)$. We give an explicit construction for $\gr(2,n)$ here. 
    
    Consider the minor whose rows are indexed by cluster variables of the form $\pl{ij}$ for $i < j-1$ and $j\geq 4$. By adding the column indexed by $\pl{jn}$ to the column $\pl{1j}$ for each $3\leq j < n-1$ the column $\pl{1j}$ becomes the column for $\gr(2,n-1)$ and thus has no entries in any row labeled by $\pl{jn}$. The remaining entries of these rows form an upper triangular matrix of the form 
    \[ \begin{blockarray}{c c  c c c c }
            & \pl{1(n-1)} & \pl{(n-2)n} & \pl{(n-3)n} & \dots & \pl{2n}\\
            \begin{block}{c [c  c c c c]}
                \pl{(n-2)n} & -1  &-1 &   & & \\
                \pl{(n-3)n} &     & 1 &-1 & & \\
                \vdots      &     &   & \ddots&\ddots  &  \\
                \pl{2n}     &     &   &  & 1 & -1\\
                \pl{1n}     &     &   &   & & 1\\
            \end{block}
        \end{blockarray} \]
    Continuing inductively this minor reduces to an upper triangular matrix with $\pm1$ on the diagonal. Each column operation leaves the determinant unchanged and so the minor has determinant $\pm1$ as needed.
\end{proof}

\begin{remark}
    The bounded ratios in $\gr(2,n)$ are both subtraction free as polynomials in all cluster variables and in the face weights of a weighted planar network. This is because every cluster variable in $\gr(k,n)$ is a positive function of the face weights.
\end{remark}

\begin{example}
    In $\gr(3,6)$ there are 16 mutable cluster variables and $6$ frozen variables. Two of the mutable variables are degree two polynomials that we call $\pl{124|356}$ and $\pl{135|246}$ using the tableaux indexing of \cite{cdfl-TableuxIndexing}.  To find the bounded Pl\"ucker cone we must find positive linear combinations of the extreme ratios (\Cref{fig:ExtremeGr36}) that eliminate these variables. If we order the Pl\"ucker variables by reverse lexicographic order and then append $\pl{124|356}, \pl{135|246}$ we obtain the following matrix for $U$:
    \begin{equation*}
        \setlength{\arraycolsep}{1pt}
        \begin{bNiceMatrix}[first-row]
           { \pl{356}}& { \pl{346}}& { \pl{256}}& { \pl{246}}& { \pl{245}}& { \pl{236}}& { \pl{235}}& { \pl{146}}& { \pl{145}}& { \pl{136}}& { \pl{135}}& { \pl{134}}& { \pl{125}}& { \pl{124}}& { \pl{124|356}}& { \pl{135|246}}\\
            1& 0& 0&  0& 0& 0& 0& 0& 0& 0& -1& 1& 1& 0& -1&  0\\
            0 &0& 0& -1& 1& 1& 0& 1& 0& 0&  0& 0& 0& 0&  0& -1\\
        \end{bNiceMatrix}
    \end{equation*}
    Clearly the ratios associated to $\pl{346},\pl{256},\pl{235},\pl{145},\pl{136}$ and $\pl{124}$ are already extreme rays of the bounded Pl\"ucker cone. We obtain 12 new extreme rays given by multiplying the extreme ray at the center of a $D_4$ with one of its three neighbors. In this way we recover the set of 18 primitive Pl\"ucker ratios.
\end{example}

\begin{proposition}
    Every Pl\"ucker bounded ratio on $\gr(3,6)$ is subtraction free. Furthermore such a ratio can be factored as a positive integer combination of primitive ratios. 
\end{proposition}
\begin{proof}
   By \Cref{thm:Gr3678SubtractionFree} we also know that every Pl\"ucker bounded ratio is subtraction free as a function of cluster variables. Since each cluster variable (including the exotic variables) is a positive function in face weights this shows \Cref{c:2} holds as well.\\
   The factorization into primitive ratios is inductive. Consider a minimal subproduct of the factorization in $\uvar{v}_\root{\gamma}$ such that the result is a Pl\"ucker ratio. By \Cref{thm:ExtremeRaysOfSBoundedCone} this is an extreme ray of the Pl\"ucker bounded cone and thus is a primitive ratio. The remaining product is smaller and so can be inductively factored.
\end{proof}

\begin{figure}[!hb]

    \begin{align*}
    \begin{aligned}
        \uvar{v}_{\pl{124|356}}=~& \frac{\pl{346} \pl{256} \pl{124}}{\pl{246}\pl{124|356}} = \dynkin D{*O**}\\
        \uvar{v}_{\pl{356}} =~& \frac{\pl{124|356}}{\pl{356} \pl{124}}= \dynkin D{O*oo}\\
        \uvar{v}_{\pl{134}} =~& \frac{\pl{124|356}}{\pl{256} \pl{134}}= \dynkin D{o*Oo}\\
        \uvar{v}_{\pl{125}} =~& \frac{\pl{124|356}}{\pl{346} \pl{125}}= \dynkin D{o*oO}
    \end{aligned}
    \hspace{2pc}
    \begin{aligned}
        v_{\pl{246}} =~&  \frac{\pl{245} \pl{236} \pl{146}}{\pl{246} \pl{135|246}} = \dynkin D{*O**}\\
        v_{\pl{124}} =~& \frac{\pl{246} \pl{123}}{\pl{236} \pl{124}} = \dynkin D{O*oo}\\
        v_{\pl{256}} =~& \frac{\pl{246} \pl{156}}{\pl{256} \pl{146}} = \dynkin D{o*Oo}\\
        v_{\pl{346}} =~& \frac{\pl{345} \pl{246}}{\pl{346} \pl{245}}= \dynkin D{o*oO}
    \end{aligned} 
 \\\\
    \begin{aligned}
        v_{\pl{135|246}} =~& \frac{\pl{235} \pl{145} \pl{136}}{\pl{135} \pl{135|246}}= \dynkin D{*O**}\\
        v_{\pl{236}} =~& \frac{\pl{135|246}}{\pl{236} \pl{145}} = \dynkin D{O*oo}\\
        v_{\pl{146}} =~& \frac{\pl{135|246}}{\pl{235} \pl{146}} = \dynkin D{o*Oo}\\
        v_{\pl{245}} =~& \frac{\pl{135|246}}{\pl{245} \pl{136}}= \dynkin D{o*oO}
    \end{aligned}
    \hspace{2pc}
    \begin{aligned}
        v_{\pl{135}} =~& \frac{\pl{356} \pl{134} \pl{125}}{\pl{135} \pl{124|356}} = \dynkin D{*O**}\\
        v_{\pl{145}} =~& \frac{\pl{456} \pl{135}}{\pl{356} \pl{145}}= \dynkin D{O*oo}\\
        v_{\pl{235}} =~& \frac{\pl{234} \pl{135}}{\pl{235} \pl{134}} = \dynkin D{o*Oo}\\
        v_{\pl{136}} =~& \frac{\pl{135} \pl{126}}{\pl{136} \pl{125}}= \dynkin D{o*oO}
    \end{aligned}
    \end{align*}
\caption{Extreme Rays in $\gr(3,6)$}
\label{fig:ExtremeGr36}
\end{figure}

\begin{example}
    In $\gr(3,7)$ there are 14 non-Pl\"ucker variables of the 42 mutable variables. Via an analogous computation to $\gr(3,6)$ we compute the cone of bounded Pl\"ucker ratios. We can classify the extreme rays of this cone into 6 orbits under the action of rotation. One representative of each orbit is factored into extreme ratios of the full cone in \Cref{fig:PluckerExtremeGr37}. As each extreme ray appears at most once, we can graphically represent the product by circling the node the Dynkin seed corresponding to each ratio. In \Cref{fig:PluckerExtremeGr37} we observe the rays of each factorization are adjacent in the sources-sink walk through the Dynkin seeds.
\end{example}

The same proof as for $\gr(3,6)$ works in $\gr(3,7)$ to show:
\begin{proposition}
    Every Pl\"ucker bounded ratio on $\gr(3,7)$ is subtraction free. Furthermore such a ratio can be factored as a positive integer combination of primitive ratios. 
\end{proposition}

\begin{figure}[!htb]
        \begin{align*}
        \begin{aligned}
            \frac{\pl{247} \pl{123}}{\pl{237} \pl{124}} =~& v_{\pl{124}}\\ 
            =~& \dynkin[labels={,,,0,,}] E{o**o*O} \\
            \frac{\pl{247} \pl{157}}{\pl{257} \pl{147}} =~& v_{\pl{257}} v_{\pl{246|357}} v_{\pl{256}}\\ 
            =~& \dynkin[labels={,,,-1,,}] E{*oO*o*} \dynkin[labels={,,,0,,}] E{o**O*o} \dynkin[labels={,,,1,,}] E{*Oo*o*} \\
            \frac{\pl{346} \pl{256}}{\pl{356} \pl{246}} =~& v_{\pl{356}} v_{\pl{124|356}}  v_{\pl{246|357}} v_{\pl{146|357}} \\
            =~&\dynkin[labels={,,,-2,,}] E{o**o*O} \dynkin[labels={,,,-1,,}] E{*oo*O*} \dynkin[labels={,,,0,,}] E{o**O*o} \dynkin[labels={,,,1,,}] E{*oO*o*} 
        \end{aligned}\hspace{1pc}
        \begin{aligned}
           \frac{\pl{256} \pl{247}}{\pl{257} \pl{246}} =~& v_{\pl{257}}( v_{\pl{157}}v_{\pl{246|357}}) v_{\pl{146|357}}\\
           =~& \dynkin[labels={,,,-1,,}] E{*oO*o*} \dynkin[labels={,,,0,,}] E{O**O*o} \dynkin[labels={,,,1,,}] E{*oO*o*} \\
           \frac{\pl{346} \pl{247}}{\pl{347} \pl{246}} =~&   v_{\pl{347}} v_{\pl{246|357}} v_{\pl{146|357}}\\
           =~&\dynkin[labels={,,,-1,,}] E{*Oo*o*} \dynkin[labels={,,,0,,}] E{o**O*o} \dynkin[labels={,,,1,,}] E{*oO*o*} \\
           \frac{\pl{256} \pl{157}}{\pl{257} \pl{156}}=~& v_{\pl{257}} v_{\pl{246|357}} v_{\pl{247}} v_{\pl{237}}\\
           =~& \dynkin[labels={,,,-1,,}] E{*oO*o*} \dynkin[labels={,,,0,,}] E{o**O*o} \dynkin[labels={,,,1,,}] E{*oo*O*} \dynkin[labels={,,,2,,}] E{o**o*O} \\
        \end{aligned}               
        \end{align*}
    \caption{Representatives of Each Orbit of Extreme Pl\"ucker Ratios in $\gr(3,7)$}
    \label{fig:PluckerExtremeGr37}
\end{figure}

\begin{example}\label{ex:Gr38PluckerCone}
    In $\gr(3,8)$ there are 80 non-Pl\"ucker variables of the 128 mutable variables. Using the computation described above we obtain 10 orbits of 8 extreme ratios under the cluster modular group. We provide the list of orbit representatives in \Cref{sec:Gr38PluckerRatios}. One orbit contains the 8 ratios in the full cone that are already Pl\"ucker. The other nine orbits consist of products of 6 or 11 extreme rays. As in $\gr(3,7)$ these products are adjacent on the sources-sink walk.  
\end{example}

Once again, the same proof as for $\gr(3,6)$ extends to $\gr(3,8)$ to show:
\begin{proposition}
    Every Pl\"ucker bounded ratio on $\gr(3,8)$ is subtraction free. Furthermore such a ratio can be factored as a positive integer combination of primitive ratios. 
\end{proposition}
\begin{example}\label{ex:Gr38WightL3Cone}
    There is a further refinement of the cluster variables in $\gr(3,8)$. The 80 non-Pl\"ucker variables can be divided into the 56 degree 2 polynomials in Pl\"ucker coordinates and the 24 degree 3 polynomials. If we let $S$ be the set of cluster variables of degree less than or equal to 2, we obtain a cone $C_S$ between the Pl\"ucker cone and the full cone of bounded ratios. Here we find 56 extreme rays of the full cone already belong to $V_S$. There are 112 new extreme rays consisting of products of 3, 4 or 5 extreme ratios of the full cone. A list of orbit representatives is given in \Cref{sec:Gr38WeightLess3Ratios}). 
\end{example}

The next biggest Grassmannian is $\gr(4,8)$. The corresponding cluster algebra is not of finite type. It was shown in  \cite{soskin2023bounded} that the Pl\"ucker bounded cone in this case is finitely generated but has extreme rays that are not primitive. In \Cref{fig:Gr48ExtremeRays} we give a full list of extreme rays up to the dihedral action and duality $I  \mapsto [8]\setminus I$. We have verified that each ratio in \Cref{fig:Gr48ExtremeRays} is subtraction-free in the sense of definition~\ref{df:subfreeFW}, which proves the following statement. 
\begin{theorem}
    \item 1. Every bounded ratio in Pl\"ucker coordinates of $\gr(4,8)$ is bounded by 1. 
    \item 2. A bounded ratio in Pl\"ucker coordinates of $\gr(4,8)$ which is a product of positive integer powers of generating ratios listed in \Cref{fig:Gr48ExtremeRays} is subtraction-free in the sense of \Cref{df:subfreeFW}.   
\end{theorem}
\begin{corollary} Let $A$ be a totally positive $4\times4$ matrix.
    If a multiplicative determinantal inequality \begin{equation*}{\det(A}_{I_{1},I'_{1}}){\det(A}_{I_{2},I'_{2}})...{\det(A}_{I_{p},I'_{p}})\leq C\cdot{\det(A}_{J_{1},J'_{1}}) {\det(A}_{J_{2},J'_{2}})...{\det(A}_{J_{q},J'_{q}})
\end{equation*}   
where $I_k,I'_k,J_k,J'_k \subseteq \{1,2,3,4\}$ with
$|I_k|=|I'_k|$ and $|J_k|=|J'_k|$ holds for some constant $C>0$ then this inequality holds with $C=1$.
\end{corollary}
\begin{figure}
    \centering
    \include{tikz/Gr48RaysPlucker}
    \caption{Extreme Pl\"ucker Rays in $\gr(4,8)$}
    \label{fig:Gr48ExtremeRays}
\end{figure}

While the sets of $U-$variables for cluster algebras of infinite type are less understood, the ratios of the form $\frac{y}{1+y}$ are naturally bounded quantities on any cluster algebra. We observe that some of these new extreme rays can be factored into the simpler ratios of the form $\frac{y}{1+y}$ for $y$ in $\gr(4,8)$. For example

{\small
\begin{equation*}
\begin{aligned}
     \frac{\pl{1368}\pl{1458}\pl{1467}\pl{2367}\pl{2358}}{\pl{1358}\pl{1367}\pl{1468}\pl{2368}\pl{2457}} =& \frac{\pl{1467}\pl{1458}\pl{1368}}{\pl{1468}\pl{1357|1468}} \frac{\pl{1357|2467}}{\pl{2457}\pl{1367}} \frac{\pl{1357|1357|2468}}{\pl{1358}\pl{1357|2467}}\\ &\frac{\pl{2367}\pl{2358}\pl{1368}}{\pl{2368}\pl{1357|2368}}\frac{\pl{1357|2368}\pl{1357|1468}}{\pl{1368}\pl{1357|1357|2468}}
\end{aligned} 
\end{equation*}
}
\section{Acknowledgements}
The authors would like to thank Thomas Lam for productive discussions and in particular for his crucial insight on the connection between bounded ratios and u-variables discussed in \cite{ahsl-ClusterConfigurationOfFiniteType}. M.G. was partially supported by NSF grant DMS $\#$2100785. He is also grateful to the Max Planck Institute for Mathematics in the Sciences for its hospitality during the June 2024 research visit. Z.G. was supported by the European Research Council under ERC-Advanced Grant 101018839.

\bibliographystyle{alpha}
\bibliography{BoundedRatioPaper}

\appendix
\renewcommand{\thesection}{\Alph{section}}
\section{Extreme Ratios in \texorpdfstring{$\gr(3,8)$}{gr 3 8}}
\subsection{Extreme \texorpdfstring{Pl\"ucker}{Plucker} Ratios in \texorpdfstring{$\gr(3,8)$}{gr 3 8}}\label{sec:Gr38PluckerRatios}
Here we give a representative of each orbit of extreme rays of the Pl\"ucker cone in $\gr(3,8)$ (\Cref{ex:Gr38PluckerCone}). For each representative we give the factorization into extreme rays of the full bounded cone ($U-$variables) and give a graphical representation of the location of the $U-$variable in the Dynkin seeds.
\input{tikz/Gr38RaysPlucker}
\subsection{Extreme Ratios in Cluster Coordinates of Degree \texorpdfstring{$\le 2$}{less than or equal to 2} in \texorpdfstring{$\gr(3,8)$}{gr 3 8}}\label{sec:Gr38WeightLess3Ratios}
Next we provide the same analysis for the cone of cluster variables of weight less than 3 in $\gr(3,8)$ (\Cref{ex:Gr38WightL3Cone}).
\input{tikz/Gr38RaysWeightLess2}

\end{document}

%% file: tikz/Gr48RaysPlucker.tex
\[
{\renewcommand{\arraystretch}{2}
\begin{array}{c}
\frac{\pl{5678} \pl{1467}}{\pl{4678} \pl{1567}} \hspace{2pc}
\frac{\pl{3578} \pl{1347}}{\pl{3478} \pl{1357}} \hspace{2pc}
\frac{\pl{3568} \pl{2578}}{\pl{3578} \pl{2568}} \hspace{2pc}
\frac{\pl{4578} \pl{3678}}{\pl{4678} \pl{3578}} \\
\frac{\pl{4578} \pl{1357}}{\pl{3578} \pl{1457}} \hspace{2pc}
\frac{\pl{3578} \pl{1356}}{\pl{3568} \pl{1357}} \hspace{2pc}
\frac{\pl{3568} \pl{3478}}{\pl{3578} \pl{3468}} \hspace{2pc}
\frac{\pl{4678} \pl{1457}}{\pl{4578} \pl{1467}}\\
\frac{\pl{3568} \pl{3478} \pl{2578} \pl{1346} \pl{1247}}{\pl{3578} \pl{3468} \pl{2478} \pl{1347} \pl{1256}} \\
\frac{\pl{3578} \pl{3468} \pl{2478} \pl{1347} \pl{1256}}{\pl{3568} \pl{3478} \pl{2468} \pl{1357} \pl{1247}} \\
\frac{\pl{3568} \pl{3478} \pl{2578} \pl{1346} \pl{1256}}{\pl{3578} \pl{3468} \pl{2568} \pl{1356} \pl{1247}} \\
\frac{\pl{4578} \pl{3678} \pl{2458} \pl{2368} \pl{2357} \pl{1467} \pl{1358}}{\pl{4678} \pl{3578} \pl{2468} \pl{2457} \pl{2358} \pl{1368} \pl{1357}} \\
\frac{\pl{4578} \pl{3678} \pl{2458} \pl{2368} \pl{2357} \pl{1467} \pl{1358}}{\pl{4678} \pl{3578} \pl{2468} \pl{2367} \pl{2358} \pl{1458} \pl{1357}} \\
\frac{\pl{4678} \pl{4568} \pl{3578} \pl{2678} \pl{2456} \pl{2357} \pl{2348} \pl{2346} \pl{1457} \pl{1367} \pl{1356} \pl{1347} \pl{1258}}{\pl{4578} \pl{3678} \pl{3568} \pl{2578} \pl{2457} \pl{2358} \pl{2356} \pl{2347} \pl{1467} \pl{1456} \pl{1348} \pl{1346} \pl{1267}} \\
\frac{\pl{4678} \pl{4568} \pl{3578} \pl{2678} \pl{2456} \pl{2357} \pl{2348} \pl{2346} \pl{1457} \pl{1367} \pl{1356} \pl{1347} \pl{1258} \pl{1258}}{\pl{4578} \pl{3678} \pl{3568} \pl{2578} \pl{2457} \pl{2358} \pl{2356} \pl{2347} \pl{1467} \pl{1456} \pl{1357} \pl{1346} \pl{1268} \pl{1248}} \\
\frac{\pl{4678} \pl{4568} \pl{3578} \pl{3468} \pl{2578} \pl{2456} \pl{2357} \pl{2348} \pl{2346} \pl{1457} \pl{1358} \pl{1356} \pl{1356} \pl{1347} \pl{1347} \pl{1267}^2 \pl{1248} \pl{1246}}{\pl{4578} \pl{3568} \pl{3478} \pl{2678} \pl{2468} \pl{2457} \pl{2358} \pl{2356} \pl{2347} \pl{1467} \pl{1456} \pl{1357}^2 \pl{1348} \pl{1346}^2 \pl{1268} \pl{1256} \pl{1247}} \\
\frac{\pl{4678} \pl{4568} \pl{3578} \pl{3468} \pl{2578} \pl{2456} \pl{2357} \pl{2348} \pl{2346} \pl{1457} \pl{1358} \pl{1356} \pl{1356} \pl{1347} \pl{1347} \pl{1267}^2 \pl{1248} \pl{1246}}{\pl{4578} \pl{3568}^2 \pl{2678} \pl{2468} \pl{2457} \pl{2358} \pl{2356} \pl{2347} \pl{1467} \pl{1456} \pl{1357}^2 \pl{1348} \pl{1346}^2 \pl{1268} \pl{1247} \pl{1247}} \\
{\frac{\pl{4678} \pl{4568} \pl{3578} \pl{3468} \pl{2578} \pl{2467} \pl{2458} \pl{2456} \pl{2357} \pl{2357} \pl{2348} \pl{2346} \pl{1457} \pl{1457} \pl{1368} \pl{1358} \pl{1356} \pl{1356} \pl{1347} \pl{1347} \pl{1267}^2 \pl{1248} \pl{1246}}{\pl{4578} \pl{3568} \pl{3478} \pl{2678} \pl{2468}^2 \pl{2457}^2 \pl{2367} \pl{2358} \pl{2356} \pl{2347} \pl{1467} \pl{1458} \pl{1456} \pl{1357}^3 \pl{1348} \pl{1346}^2 \pl{1268} \pl{1256} \pl{1247}} }\\
\frac{\pl{4678} \pl{4568} \pl{3578} \pl{3468} \pl{2578} \pl{2467} \pl{2458} \pl{2456} \pl{2357} \pl{2357} \pl{2348} \pl{2346} \pl{1457} \pl{1457} \pl{1368} \pl{1358} \pl{1356} \pl{1356} \pl{1347} \pl{1347} \pl{1267}^2 \pl{1248} \pl{1246}}{\pl{4578} \pl{3568}^2 \pl{2678} \pl{2468}^2 \pl{2457}^2 \pl{2358}^2 \pl{2356} \pl{2347} \pl{1467}^2 \pl{1456} \pl{1357}^3 \pl{1348} \pl{1346}^2 \pl{1268} \pl{1247}^2} \\
\frac{\pl{4678} \pl{4568} \pl{3578} \pl{3468} \pl{2578} \pl{2467} \pl{2458} \pl{2456} \pl{2357} \pl{2357} \pl{2348} \pl{2346} \pl{1457}^2 \pl{1368} \pl{1358} \pl{1356}^2 \pl{1347}^2 \pl{1267}^2 \pl{1248} \pl{1246}}{\pl{4578} \pl{3568} \pl{3478} \pl{2678} \pl{2468}^2 \pl{2457}^2 \pl{2358}^2 \pl{2356} \pl{2347} \pl{1467}^2 \pl{1456} \pl{1357}^3 \pl{1348} \pl{1346}^2 \pl{1268} \pl{1256} \pl{1247}} \\
\end{array}
}
\]

%% file: tikz/Gr38RaysPlucker.tex
\begingroup
\allowdisplaybreaks
\begin{align*}
        \frac{\pl{456}\pl{357}}{\pl{457}\pl{356}} =~& v_{\pl{457}}
         =~ \dynkin[labels={,,,0,,,,}] E{o**o*o*O}\\
        \frac{\pl{258} \pl{168}}{\pl{268} \pl{158}}=~& v_{\pl{268}}v_{\pl{257|468}} v_{\pl{257|368}} v_{\pl{124|257|368}}v_{\pl{247|368}} v_{\pl{267}}\\
        =~& \dynkin[labels={,,,-3,,,,}] E{*oo*o*O*} \dynkin[labels={,,,-2,,,,}] E{o**o*O*o} \dynkin[labels={,,,-1,,,,}] E{*oo*O*o*} \dynkin[labels={,,,0,,,,}] E{o**O*o*o} \dynkin[labels={,,,1,,,,}] E{*oO*o*o*}\dynkin[labels={,,,2,,,,}] E{O**o*o*o}\\
        \frac{\pl{258} \pl{124}}{\pl{248} \pl{125}} =~&v_{\pl{125}} v_{\pl{124|358}}v_{\pl{124|257|368}} v_{\pl{124|157|368}}v_{\pl{124|357}}v_{\pl{124|356}}\\
        =~&  \dynkin[labels={,,,-2,,,,}] E{O**o*o*o} \dynkin[labels={,,,-1,,,,}] E{*oO*o*o*} \dynkin[labels={,,,0,,,,}] E{o**O*o*o} \dynkin[labels={,,,1,,,,}] E{*oo*O*o*}\dynkin[labels={,,,2,,,,}] E{o**o*O*o}          \dynkin[labels={,,,3,,,,}] E{*oo*o*O*}\\
        \frac{\pl{348} \pl{124}}{\pl{248} \pl{134}} =~&v_{\pl{134}}v_{\pl{124|378}} v_{\pl{124|368}}v_{\pl{124|357|468}}v_{\pl{124|357|368}}\\
          &(v_{\pl{124|358}}v_{\pl{124|367}})v_{\pl{124|257|368}}v_{\pl{124|157|368}}v_{\pl{124|357}} v_{\pl{124|356}}\\
           =~& \dynkin[labels={,,,-6,,,,}] E{o**o*o*O} \dynkin[labels={,,,-5,,,,}] E{*oo*o*O*}\dynkin[labels={,,,-4,,,,}] E{o**o*O*o} \dynkin[labels={,,,-3,,,,}] E{*oo*O*o*}\dynkin[labels={,,,-2,,,,}] E{o**O*o*o}\\
           &\dynkin[labels={,,,-1,,,,}]E{*OO*o*o*} \dynkin[labels={,,,0,,,,}] E{o**O*o*o} \dynkin[labels={,,,1,,,,}] E{*oo*O*o*}\dynkin[labels={,,,2,,,,}] E{o**o*O*o} \dynkin[labels={,,,3,,,,}] E{*oo*o*o*}\\
        \frac{\pl{267} \pl{168}}{\pl{268} \pl{167}} =~& v_{\pl{268}} v_{\pl{257|468}} v_{\pl{257|368}}v_{\pl{124|257|368}}  (v_{\pl{247|368}}v_{\pl{258}})\\
        & v_{\pl{147|258|368}}v_{\pl{247|358}}v_{\pl{246|358}}v_{\pl{248}}v_{\pl{238}}\\ 
      =~&  \dynkin[labels={,,,-3,,,,}] E{*oo*o*O*}\dynkin[labels={,,,-2,,,,}] E{o**o*O*o} \dynkin[labels={,,,-1,,,,}] E{*oo*O*o*}\dynkin[labels={,,,0,,,,}] E{o**O*o*o} \dynkin[labels={,,,1,,,,}] E{*OO*o*o*}\\& \dynkin[labels={,,,2,,,,}] E{o**O*o*o} \dynkin[labels={,,,3,,,,}] E{*oo*O*o*}\dynkin[labels={,,,4,,,,}] E{o**o*O*o} \dynkin[labels={,,,5,,,,}] E{*oo*o*O*}\dynkin[labels={,,,6,,,,}] E{o**o*o*O}\\
              \frac{\pl{357} \pl{267}}{\pl{367} \pl{257}} =~& v_{\pl{367}} v_{\pl{125|367}} v_{\pl{124|357|368}} (v_{ \pl{257|368}} v_{ \pl{124|367}}) \\
              & (v_{\pl{124|257|368}}v_{\pl{157|368}})(v_{ \pl{247|368}} v_{\pl{124|157|368}}) v_{\pl{147|258|368}} v_{\pl{147|368}}\\
            =~& \dynkin[labels={,,,-4,,,,}] E{O**o*o*o}\dynkin[labels={,,,-3,,,,}] E{*oO*o*o*}\dynkin[labels={,,,-2,,,,}] E{o**O*o*o}\dynkin[labels={,,,-1,,,,}] E{*Oo*O*o*} \\& \dynkin[labels={,,,0,,,,}] E{o**O*O*o}\dynkin[labels={,,,1,,,,}] E{*oO*O*o*}\dynkin[labels={,,,2,,,,}] E{o**O*o*o}\dynkin[labels={,,,3,,,,}] E{*Oo*o*o*}\\
        \frac{\pl{357} \pl{348}}{\pl{358} \pl{347}} =~&  v_{\pl{358}}v_{\pl{124|357|368}}(v_{\pl{124|358}}v_{\pl{257|368}})(v_{\pl{124|257|368}}v_{\pl{157|368}})\\
            &(v_{\pl{124|157|368}}v_{\pl{258}}) v_{\pl{147|258|368}}v_{\pl{157|268}}v_{\pl{158}}  \\
            =~&\dynkin[labels={,,,-3,,,,}] E{*Oo*o*o*}\dynkin[labels={,,,-2,,,,}] E{o**O*o*o}\dynkin[labels={,,,-1,,,,}] E{*oO*O*o*}\dynkin[labels={,,,0,,,,}] E{o**O*O*o} \\&\dynkin[labels={,,,1,,,,}] E{*Oo*O*o*}\dynkin[labels={,,,2,,,,}] E{o**O*o*o}\dynkin[labels={,,,3,,,,}] E{*oO*o*o*}\dynkin[labels={,,,4,,,,}] E{O**o*o*o}\\
        \frac{\pl{348} \pl{258}}{\pl{358} \pl{248}} =~& v_{\pl{358}}v_{\pl{124|357|368}}(v_{\pl{124|358}}v_{\pl{257|368}})(v_{\pl{124|257|368}}v_{\pl{157|368}})\\
            &(v_{\pl{124|157|368}}v_{\pl{357}})(v_{\pl{124|357}}v_{\pl{356}})v_{\pl{124|356}}\\
            =~& \dynkin[labels={,,,-3,,,,}] E{*Oo*o*o*}\dynkin[labels={,,,-2,,,,}] E{o**O*o*o}\dynkin[labels={,,,-1,,,,}] E{*oO*O*o*}\dynkin[labels={,,,0,,,,}] E{o**O*O*o} \\&\dynkin[labels={,,,1,,,,}] E{*oo*O*O*}\dynkin[labels={,,,2,,,,}] E{o**o*O*O}\dynkin[labels={,,,3,,,,}] E{*oo*o*O*}\\
        \frac{\pl{267} \pl{258}}{\pl{268} \pl{257}} =~& v_{\pl{268}}(v_{\pl{257|468}}v_{\pl{168}})(v_{\pl{257|368}}  v_{\pl{157|468}})(v_{\pl{124|257|368}}v_{\pl{157|368}})\\
            &(v_{\pl{247|368}}v_{\pl{124|157|368}})v_{\pl{147|258|368}} v_{\pl{147|368}}     \\ 
            =~& \dynkin[labels={,,,-3,,,,}] E{*oo*o*O*}\dynkin[labels={,,,-2,,,,}] E{o**o*O*O}\dynkin[labels={,,,-1,,,,}] E{*oo*O*O*}\dynkin[labels={,,,0,,,,}] E{o**O*O*o} \\&\dynkin[labels={,,,1,,,,}] E{*oO*O*o*}\dynkin[labels={,,,2,,,,}] E{o**O*o*o}\dynkin[labels={,,,3,,,,}] E{*Oo*o*o*}\\
        \frac{\pl{357} \pl{258}}{\pl{358} \pl{257}} =~& v_{\pl{358}}v_{\pl{124|357|368}}(v_{\pl{124|358}}v_{\pl{257|368}})(v_{\pl{348}}v_{\pl{124|257|368}}v_{\pl{157|368}})\\
        &(v_{\pl{247|368}}v_{\pl{124|157|368}})v_{\pl{147|258|368}}v_{\pl{147|368}} \\
        =~& \dynkin[labels={,,,-3,,,,}] E{*Oo*o*o*}\dynkin[labels={,,,-2,,,,}] E{o**O*o*o}\dynkin[labels={,,,-1,,,,}] E{*oO*O*o*}\dynkin[labels={,,,0,,,,}] E{O**O*O*o} \\&\dynkin[labels={,,,1,,,,}] E{*oO*O*o*}\dynkin[labels={,,,2,,,,}] E{o**O*o*o}\dynkin[labels={,,,3,,,,}] E{*Oo*o*o*}
\end{align*}
\endgroup

%% file: tikz/Gr38RaysWeightLess2.tex
\begingroup
\allowdisplaybreaks
\begin{align*}
      \frac{\pl{247|358} \pl{157|268}}{\pl{258} \pl{257} \pl{147|368}} =~&v_{\pl{258}} v_{\pl{147|258|368}} v_{\pl{147|368}}
      = \dynkin[labels={,,,1,,,,}] E{*Oo*o*o*} \dynkin[labels={,,,2,,,,}] E{o**O*o*o} \dynkin[labels={,,,3,,,,}] E{*Oo*o*o*}\\
      \frac{\pl{158} \pl{247|358}}{\pl{258} \pl{147|358}} =~&v_{\pl{258}} v_{\pl{147|258|368}} v_{\pl{157|268}}
      = \dynkin[labels={,,,1,,,,}] E{*Oo*o*o*} \dynkin[labels={,,,2,,,,}] E{o**O*o*o} \dynkin[labels={,,,3,,,,}] E{*oO*o*o*}\\
      \frac{\pl{267} \pl{247|358}}{\pl{257} \pl{247|368}} =~&v_{\pl{247|368}} v_{\pl{147|258|368}} v_{\pl{147|368}}
      = \dynkin[labels={,,,1,,,,}] E{*oO*o*o*} \dynkin[labels={,,,2,,,,}] E{o**O*o*o} \dynkin[labels={,,,3,,,,}] E{*Oo*o*o*}\\
      \frac{\pl{247|368} \pl{157|368}}{\pl{257|368} \pl{147|368}} =~&v_{\pl{257|368}} v_{\pl{124|257|368}} v_{\pl{258}}
      = \dynkin[labels={,,,-1,,,,}] E{*oo*O*o*}\dynkin[labels={,,,0,,,,}] E{o**O*o*o}  \dynkin[labels={,,,1,,,,}] E{*Oo*o*o*}\\ 
      \frac{\pl{247|368} \pl{124|357}}{\pl{247|358} \pl{124|367}} =~&v_{\pl{124|367}} v_{\pl{124|257|368}} v_{\pl{124|157|368}}
      = \dynkin[labels={,,,-1,,,,}] E{*Oo*o*o*}\dynkin[labels={,,,0,,,,}] E{o**O*o*o}  \dynkin[labels={,,,1,,,,}] E{*Oo*o*o*}\\
      \frac{\pl{248} \pl{157|268}}{\pl{258} \pl{147|268}} =~&v_{\pl{258}} v_{\pl{147|258|368}} v_{\pl{247|358}} v_{\pl{246|358}}\\
      =& \dynkin[labels={,,,1,,,,}] E{*Oo*o*o*} \dynkin[labels={,,,2,,,,}] E{o**O*o*o} \dynkin[labels={,,,3,,,,}] E{*oo*O*o*}\dynkin[labels={,,,4,,,,}] E{o**o*O*o}\\
      \frac{\pl{357} \pl{157|268}}{\pl{257} \pl{157|368}} =~&v_{\pl{157|368}} v_{\pl{124|157|368}} v_{\pl{147|258|368}} v_{\pl{147|368}}\\
      =&\dynkin[labels={,,,0,,,,}] E{o**o*O*o}  \dynkin[labels={,,,1,,,,}] E{*oo*O*o*} \dynkin[labels={,,,2,,,,}] E{o**O*o*o} \dynkin[labels={,,,3,,,,}] E{*Oo*o*o*}\\
      \frac{\pl{267} \pl{158} \pl{247|358} \pl{147|368}}{\pl{247|368} \pl{147|358} \pl{157|268}} =~&v_{\pl{247|368}} v_{\pl{147|258|368}} v_{\pl{157|268}}
      = \dynkin[labels={,,,1,,,,}] E{*oO*o*o*} \dynkin[labels={,,,2,,,,}] E{o**O*o*o} \dynkin[labels={,,,3,,,,}] E{*oO*o*o*} \\
      \frac{\pl{267} \pl{258} \pl{157|368}}{\pl{257|368} \pl{157|268}} =~&v_{\pl{257|368}} v_{\pl{124|257|368}} v_{\pl{247|368}}
      = \dynkin[labels={,,,-1,,,,}] E{*oo*O*o*}\dynkin[labels={,,,0,,,,}] E{o**O*o*o}  \dynkin[labels={,,,1,,,,}] E{*oO*o*o*}\\
      \frac{\pl{348} \pl{258} \pl{124|357}}{\pl{247|358} \pl{124|358}} =~&v_{\pl{124|358}} v_{\pl{124|257|368}} v_{\pl{124|157|368}}
      =\dynkin[labels={,,,-1,,,,}] E{*oO*o*o*}\dynkin[labels={,,,0,,,,}] E{o**O*o*o}  \dynkin[labels={,,,1,,,,}] E{*oo*O*o*}\\
      \frac{\pl{267} \pl{248} \pl{147|368}}{\pl{247|368} \pl{147|268}} =~&v_{\pl{247|368}} v_{\pl{147|258|368}} v_{\pl{247|358}} v_{\pl{246|358}}\\
      =& \dynkin[labels={,,,1,,,,}] E{*oO*o*o*} \dynkin[labels={,,,2,,,,}] E{o**O*o*o} \dynkin[labels={,,,3,,,,}] E{*oo*O*o*}\dynkin[labels={,,,4,,,,}] E{o**o*O*o}\\
      \frac{\pl{357} \pl{158} \pl{147|368}}{\pl{157|368} \pl{147|358}} =~&v_{\pl{157|368}} v_{\pl{124|157|368}} v_{\pl{147|258|368}} v_{\pl{157|268}}\\
      =&\dynkin[labels={,,,0,,,,}] E{o**o*O*o}  \dynkin[labels={,,,1,,,,}] E{*oo*O*o*} \dynkin[labels={,,,2,,,,}] E{o**O*o*o} \dynkin[labels={,,,3,,,,}] E{*oO*o*o*}\\
      \frac{\pl{357} \pl{258} \pl{247|368}}{\pl{257|368} \pl{247|358}} =~&v_{\pl{257|368}} (v_{\pl{157|368}} v_{\pl{124|257|368}}) v_{\pl{124|157|368}}\\
      =&\dynkin[labels={,,,-1,,,,}] E{*oo*O*o*}\dynkin[labels={,,,0,,,,}] E{o**O*O*o}  \dynkin[labels={,,,1,,,,}] E{*oo*O*o*}\\
      \frac{\pl{357} \pl{248} \pl{147|368} \pl{157|268}}{\pl{247|358} \pl{157|368} \pl{147|268}} =~&v_{\pl{157|368}} v_{\pl{124|157|368}} v_{\pl{147|258|368}} v_{\pl{247|358}} v_{\pl{246|358}}\\
      =&\dynkin[labels={,,,0,,,,}] E{o**o*O*o}  \dynkin[labels={,,,1,,,,}] E{*oo*O*o*} \dynkin[labels={,,,2,,,,}] E{o**O*o*o} \dynkin[labels={,,,3,,,,}] E{*oo*O*o*}\dynkin[labels={,,,4,,,,}] E{o**o*O*o} 
    \end{align*}
\endgroup